\documentclass[11pt]{amsart}

\usepackage{amssymb,amsthm,amsmath}
\usepackage[numbers,sort&compress]{natbib}
\usepackage{color}
\usepackage{graphicx}
\usepackage{tikz}
\usepackage[colorlinks,
            linkcolor=blue,
            anchorcolor=blue,
            citecolor=blue
            ]{hyperref}

\newcommand{\SL}{\mathrm{SL}(2,\mathbb{R})}

\newcommand{\rot}{\mathrm{rot}}

\DeclareMathOperator{\im}{Im}

\def\a{\alpha}

\def\e{\varepsilon}

\let\newpf\proof \let\proof\relax 
\newenvironment{pf}{\newpf[\proofname]}{\qed\endtrivlist}

\newcommand{\ba}{\overline{A}}

\def\be{\begin{equation}}
\def\ee{\end{equation}}

\def\ba{{\begin{align}}}
\def\ea{{\end{align}}}

\def\bm{\begin{matrix}}
\def\em{\end{matrix}}

\renewcommand{\sl}{{\mathrm{sl}}}

\def\a{{\alpha}}

\def\SL{{\mathrm{SL}}}

\def\0{{\mathbf 0}}

\newtheorem{Theorem}{Theorem}[section]
\newtheorem*{Theorem*}{Theorem 1.1}
\newtheorem{Lemma}{Lemma}[section]
\newtheorem{Proposition}{Proposition}[section]
\newtheorem{Corollary}{Corollary}[section]

\newtheorem{Definition}{Definition}[section]
\newtheorem{Claim}{Claim}

\numberwithin{equation}{section}

\theoremstyle{definition}

\newtheorem{Remark}{Remark}[section]

\def\ssm{\smallsetminus}
\def\tr{{\text{tr}}}
\renewcommand{\setminus}{\ssm}

\renewcommand{\mod}{\operatorname{mod}}

\newcommand{\D}{{\mathbb D}}

\newcommand{\R}{{\mathbb R}}
\newcommand{\T}{{\mathbb T}}

\newcommand{\X}{{\mathbb X}}
\newcommand{\Z}{{\mathbb Z}}

\def\B0{{\bold{0}}}

\catcode`\@=12

\def\Empty{}
\newcommand\oplabel[1]{
  \def\OpArg{#1} \ifx \OpArg\Empty {} \else
    \label{#1}
  \fi}

\newcommand{\comm}[1]{}
\newcommand{\comment}[1]{}

\begin{document}

\title[Cantor Spectrum for Some Almost Periodic CMV Matrices]{Cantor Spectrum for CMV Matrices With Almost Periodic Verblunsky Coefficients}

\author{Long Li}
\address{Department of Mathematics, Nanjing University, Nanjing 210093, China}
\email{huanzhensu@icloud.com}

\author{David Damanik}
\address{ Department of Mathematics, Rice University, Houston, Texas, 77005}
\email{damanik@rice.edu}
\thanks{D.D.\ was supported in part by NSF grants DMS--1700131 and DMS--2054752, an Alexander von Humboldt Foundation research award, and Simons Fellowship $\# 669836$}

\author{Qi Zhou}
\address{
Chern Institute of Mathematics and LPMC, Nankai University, Tianjin 300071, China
}
\email{qizhou@nankai.edu.cn}
\thanks{Q.Z.\ was  supported by National Key R\&D Program of China (2020YFA0713300), NSFC grant (12071232), The Science Fund for Distinguished Young Scholars of Tianjin (No. 19JCJQJC61300) and Nankai Zhide Foundation.}

\setcounter{tocdepth}{1}



\begin{abstract}
We consider extended CMV matrices with analytic quasi-periodic Verblunsky coefficients with Diophantine frequency vector in the perturbatively small coupling constant regime and prove the analyticity of the tongue boundaries. As a consequence we establish that, generically, all gaps of the spectrum that are allowed by the Gap Labelling Theorem are open and hence the spectrum is a Cantor set. We also prove these results for a related class of almost periodic Verblunsky coefficients and present an application to suitable quantum walk models on the integer lattice.
\end{abstract}

\maketitle

\section{Introduction}

CMV matrices are canonical matrix representations of unitary operators with a cyclic vector, and they arise naturally in the context of orthogonal polynomials on the unit circle. We refer the reader to \cite{Simon1, Simon2} for background. The spectral analysis of CMV matrix has seen exciting progress in the past two decades, partly due to the connection with the well-developed theory of one-dimensional discrete Schr\"odinger operators, and more generally Jacobi matrices. This exceedingly useful correspondence forms the basis of much of \cite{Simon2}.

This paper is concerned with the structure of the spectrum of CMV matrices with quasi-periodic (and related almost periodic) Verblunsky coefficients. In an earlier companion paper, \cite{LDZ}, we studied the spectral type of these operators and confirmed the expectation that it is purely absolutely continuous under suitable assumptions. Here we address the equally natural expectation that the gaps of the spectrum should generically be dense; indeed, we show that generically, all gaps allowed by the Gap Labelling Theorem are open, which in turn implies that the gaps are dense since the labels are dense and the rotation number strictly increases on the spectrum.

\subsection{Background}

Let us recall how CMV matrices arise in connection with orthogonal polynomials on the unit circle. Suppose $\nu$ is a non-trivial probability measure on the unit circle $\partial \mathbb{D} = \{ z \in \mathbb{C} : |z| = 1 \}$, that is, $\nu(\partial \D) = 1$ and $\nu$ is not supported on a finite set. By the non-triviality assumption, the functions $1$, $z$, $z^2,\cdots$ are linearly independent in the Hilbert space $\mathcal{H} = L^2(\partial\mathbb{D}, d\nu)$, and hence one can form, by the Gram-Schmidt procedure, the \textit{monic orthogonal polynomials} $\Phi_n(z)$, whose \textit{Szeg\H{o} dual} is defined by $\Phi_n^{*} = z^n\overline{\Phi_n({1}/{\overline{z}})}$. There are complex numbers $\{\alpha_n\}_{n=0}^\infty$ in $\mathbb{D}=\{z\in\mathbb{C}:|z|<1\}$, called the \textit{Verblunsky coefficients}, so that
\begin{equation}\label{eq1}
\Phi_{n+1}(z) = z \Phi_n(z) - \overline{\alpha}_n \Phi_n^*(z),
\end{equation}
which is the so-called \textit{Szeg\H{o} recurrence}. Conversely, every sequence $\{\alpha_n\}_{n=0}^\infty$ in $\mathbb{D}$ arises in this way.

The orthogonal polynomials may or may not form a basis of $\mathcal{H}$. However, if we apply the Gram-Schmidt procedure to $1, z, z^{-1}, z^2, z^{-2}, \ldots$, we will obtain a basis -- called the \textit{CMV basis}. In this basis, multiplication by the independent variable $z$ in $\mathcal{H}$ has the matrix representation
\begin{equation*}
\mathcal{C}=\left(
\begin{matrix}
\overline{\alpha}_0&\overline{\alpha}_1\rho_{0}&\rho_1\rho_0&0&0&\cdots&\\
\rho_0&-\overline{\alpha}_1\alpha_{0}&-\rho_1\alpha_0&0&0&\cdots&\\
0&\overline{\alpha}_2\rho_{1}&-\overline{\alpha}_2\alpha_{1}&\overline{\alpha}_3\rho_2&\rho_3\rho_2&\cdots&\\
0&\rho_2\rho_{1}&-\rho_2\alpha_{1}&-\overline{\alpha}_3\alpha_2&-\rho_3\alpha_2&\cdots&\\
0&0&0&\overline{\alpha}_4\rho_3&-\overline{\alpha}_4\alpha_3&\cdots&\\
\cdots&\cdots&\cdots&\cdots&\cdots&\cdots&
\end{matrix}
\right),
\end{equation*}
where
\begin{equation}\label{e.rhodef}
\rho_n = (1-|\alpha_n|^2)^{1/2}
\end{equation}
for $n \ge 0$. A matrix of this form is called a \textit{CMV matrix}.

It is sometimes helpful to also consider a two-sided extension of a matrix of this form. Namely, given a bi-infinite sequence $\{ \alpha_n \}_{n \in \Z}$ in $\D$ (and defining the $\rho_n$'s as in \eqref{e.rhodef}), we may consider the \textit{extended CMV matrix}
\begin{equation*}
\mathcal{E}=\left(
\begin{matrix}
\cdots&\cdots&\cdots&\cdots&\cdots&\cdots&\cdots\\
\cdots&-\overline{\alpha}_0\alpha_{-1}&\overline{\alpha}_1\rho_{0}&\rho_1\rho_0&0&0&\cdots&\\
\cdots&-\rho_0\alpha_{-1}&-\overline{\alpha}_1\alpha_{0}&-\rho_1\alpha_0&0&0&\cdots&\\
\cdots&0&\overline{\alpha}_2\rho_{1}&-\overline{\alpha}_2\alpha_{1}&\overline{\alpha}_3\rho_2&\rho_3\rho_2&\cdots&\\
\cdots&0&\rho_2\rho_{1}&-\rho_2\alpha_{1}&-\overline{\alpha}_3\alpha_2&-\rho_3\alpha_2&\cdots&\\
\cdots&0&0&0&\overline{\alpha}_4\rho_3&-\overline{\alpha}_4\alpha_3&\cdots&\\
\cdots&\cdots&\cdots&\cdots&\cdots&\cdots&\cdots
\end{matrix}\right).
\end{equation*}
This is especially natural in cases where the Verblunsky coefficients are generated by an invertible ergodic transformation. In this setting there is a natural extension of the definition of $\alpha_n$ to $n \in \Z_-$ and, more importantly, the spectral theory of the extended case is more elegant. In the present paper, this ``elegance'' will mean that we can obtain the Cantor structure of the spectrum of the extended CMV matrices, while the spectra of their half-line cousins will in general not be Cantor sets.\footnote{On the other hand, it is known that in our setting, the spectrum of the extended CMV matrix coincides with the essential spectrum of the standard CMV matrix, and hence the failure of Cantor spectrum for the latter is solely due to the presence of isolated eigenvalues in gaps of the essential spectrum.}

To motivate our study of the case of quasi-periodic and almost periodic Verblunsky coefficients, let us recall that in the Remarks and Historical Notes to \cite[Section~10.16]{Simon2}, Simon writes that from his discussion of ergodic Verblunsky coefficients ``conspicuously absent is the case of almost periodic Verblunsky coefficients'' and ``especially interesting is the quasiperiodic case.'' Based on this remark, there has been quite a bit of activity devoted to this issue (compare, e.g., \cite{FDG, F17, FO17, J, LDZ, O12, O13, WD, Zhang} and references therein), and in the present paper we address one of the natural problems one seeks to understand in this context: the generic presence of all spectral gaps under suitable assumptions.

\subsection{Generic Gap Opening in the Quasi-Periodic Setting}

Let us now describe the quasi-periodic setting and one of our main results in detail.

For $\rho > 0$ and a bounded analytic (possibly matrix-valued) function $F(x)$ defined on $ \{ x :  | \Im x | < \rho \}$, let $\|F\|_\rho =  \sup_{ | \Im x |< \rho } \| F(x)\| $. We denote by $C^\omega_{\rho}(\T^d,*)$ the space of all $*$-valued functions of this kind ($*$ will usually denote $\R$, $\mathrm{SU}(1,1)$, $\mathrm{SL}(2,\R)$). We also write $C^\omega(\T^d,*) = \bigcup_{\rho > 0} C^\omega_{\rho}(\T^d,*)$.

We consider analytic quasi-periodic Verblunsky coefficients of the form
\begin{equation}\label{setting}
\alpha_{n}(x) = f(x+n\omega), \quad n\in \mathbb{Z},
\end{equation}
where
\begin{equation}\label{setting2}
f(x) = \lambda e^{i h(x)},
\end{equation}
$h \in C^{\omega}(\mathbb{T}^{d},\mathbb{R})$, $\lambda \in (0,1)$, $\omega \in\mathbb{R}^{d}$ with $\langle m,\omega\rangle\notin \mathbb{Z}$ for any $m\in \mathbb{Z}^{d}\backslash\{0\}$.\footnote{Note that the symbol $\omega$ denotes both to the translation vector and the analyticity of function $h$. This double usage is standard and we hope it won't lead to any confusion.} 

The sequence of Verblunsky coefficients defined by \eqref{setting}--\eqref{setting2} determines an extended CMV matrix $\mathcal{E}$. The spectrum of $\mathcal{E}$ is independent of $x$ due to the minimality of the translation by $\omega$ on $\mathbb{T}^{d}$. It of course depends on $h$, $\lambda$, $\omega,$ and we will denote it by $\Sigma = \Sigma_{h, \lambda, \omega}$, depending on whether we want to make the parameters explicit or not. More generally, we will make only those parameters explicit in our notation that are not fixed during the discussion in question.


As $\Sigma \subseteq \partial \mathbb{D}$, its complement $\partial \mathbb{D} \setminus \Sigma$ is an at most countable union of open arcs, commonly referred to as \textit{gaps} of $\Sigma$. It is well known (cf., e.g., \cite{DFLY}) that the gaps are given by those values of the spectral parameter for which the associated Szeg\H{o} cocycle is uniformly hyperbolic. Moreover, the rotation number of this cocycle is constant in each gap and hence the value the rotation number takes in a given gap may be used to label it. The Gap Labelling Theorem (cf., e.g., \cite{GJ96, JM82}) states that the possible values the rotation number may take in gaps is fully determined by the base dynamics, which in our case is the set
$$
\mathrm{Labels}(\omega) = \left\{ \frac{\langle k,\omega\rangle}{2} \mod \mathbb{Z} : k \in  \mathbb{Z}^d \right\}.
$$
It is very natural to ask whether a possible gap label actually occurs as the value of the rotation number in a gap of $\Sigma$. Keeping $\omega$ and $\lambda$ fixed, one has an $h$-dependent spectrum $\Sigma = \Sigma_h$ and can consider the set
$$
\mathcal{O}(k) = \left\{ h \in C^{\omega}(\mathbb{T}^{d},\mathbb{R}) : \Sigma_h \text{ has a gap with rotation number } \frac{\langle k,\omega\rangle}{2} \mod \mathbb{Z} \right\}
$$
for $k \in \mathbb{Z}^d$. It is easy to see that $\mathcal{O}(k)$ is open. The much harder question is whether it is dense, or at least dense in some portion of $C^{\omega}(\mathbb{T}^{d},\mathbb{R})$. Should that indeed be the case for every $k \in \mathbb{Z}^d$, then taking the intersection over $k$ yields a $G_\delta$ set that is dense (or at least dense in the portion of $C^{\omega}(\mathbb{T}^{d},\mathbb{R})$ in question), and for each $h$ in this intersection, one can then conclude that all possible gap labels occur. Again, since they are dense, the gaps must be dense, and hence this is a stronger version of the Cantor spectrum phenomenon. In common parlance, one says that ``generically, all gaps are open.'' Our goal is to establish results of this kind.

Recall that $\omega\in\mathbb{R}^{d}$ is called {\it Diophantine} if there exist some $\kappa, \tau > 0$ such that
\begin{equation}\label{eq.diophantine}
\inf \limits_{j \in\mathbb{Z}} \vert \langle n,\omega\rangle - j \vert \geq \frac{\kappa}{\vert n\vert^{\tau}}
\end{equation}
for all $n \in \mathbb{Z}^{d} \backslash \{0\}$. Let $\mathrm{DC}(\kappa,\tau)$ be the set of all Diophantine $\omega \in \mathbb{R}^d$ with prescribed $\kappa,\tau$.

Then we have the following result, which shows that all gaps are open generically for perturbatively small $h$:

\begin{Theorem}\label{MainThm}
Let $\omega\in \mathrm{DC}(\kappa,\tau)$, $\lambda \in (0,1)$, and $ \rho>0$ be fixed. Then there exists a constant $\epsilon=\epsilon(\kappa,\tau,\rho) > 0$ such that generically in
$$
\{ h \in C^{\omega}_{\rho}(\mathbb{T}^{d},\mathbb{R}) : \|h\|_{\rho}\leq \epsilon\},
$$
all gaps of the set $\Sigma_h$ associated with the Verblunsky coefficients given by \eqref{setting}--\eqref{setting2} are open and hence for generic $h$'s in this set, $\Sigma_h$ is a Cantor set.
\end{Theorem}

\begin{Remark}
{\rm (a)} As usual, a statement holds generically if it is true for every element from a dense $G_\delta$ subset. In the present setting the dense $G_\delta$ subset will arise via the intersection over the countable set of gap labels and the key observation that for each fixed label, the corresponding gap is open for $h$'s from an open and dense set.

{\rm (b)} It is of interest to prove results of this kind in other regularity classes. In the continuous category (i.e., for $ h \in C(\mathbb{T}^{d},\mathbb{R})$), generic Cantor spectrum was shown by Jun \cite{J}.

{\rm (c)} As pointed out above, much of the recent progress in the study of orthogonal polynomials on the unit circle and CMV matrices is made possible by the analogy to the spectral analysis of one-dimensional Schr\"odinger operators. An analog of Theorem~\ref{MainThm} in the latter setting was shown by Eliasson \cite{Eli92} and Puig-Sim\'o \cite{Puig06}. We emphasize that working out the CMV analog of a Schr\"odinger result is in many cases not straightforward and one needs to overcome a number of obstacles in the implementation of the general proof strategy, which however does often resemble that from the Schr\"odinger case. This is precisely the situation we find ourselves in: we follow the general proof strategy from \cite{Puig06}, but need to deal with the resulting difficulties when implementing it.

{\rm (d)} In fact, Jun's result \cite{J} mentioned in part (b) of the present remark can be viewed as the CMV analog of a result shown in the Schr\"odinger case by Avila-Bochi-Damanik \cite{ABD09}. In this connection it is worth pointing out that the method from \cite{ABD09, J} only establishes Cantor spectrum, not the stronger version where all gaps are open. The stronger version was obtained in the Schr\"odinger setting by Avila-Bochi-Damanik in their follow-up paper \cite{ABD12}, but no CMV analog of this work is available yet. It would be of interest to work out a CMV analog of \cite{ABD12}.

{\rm (e)} Historically, the most well-known ``all gaps are open'' result is  known as the ``Dry Ten Martini Problem.'' In 1981, during a talk at the   AMS annual meeting, Mark Kac asked famously whether for the almost Mathieu operator
\begin{equation*}
(H_{\lambda,\alpha,\theta} u)_n= u_{n+1}+u_{n-1} +2\lambda \cos(
n\alpha + \theta) u_n,
\end{equation*}
``all gaps are there'' \cite{Kac,sim}, and offered ten Martinis for the solution.  Barry Simon then popularized this question in the form of two problems \cite{sim}.  The first one, called the Ten Martini Problem, asks whether the spectrum is a Cantor set. The second more difficult one, called the Dry Ten Martini Problem, asks more precisely whether all gaps allowed by the Gap Labelling Theorem are open.  While the Ten Martini Problem was completely solved by Avila-Jitomirskaya \cite{AJ05}, the Dry Ten Martini Problem remains open in the critical case $\lambda=1$ \cite{ayz}. The interested reader can consult \cite{AJ05,ayz} for more on the history of these two problems.
\end{Remark}

\subsection{Generic Gap Opening in the Almost-Periodic Setting and an Application to Quantum  Walks}

We are also interested in a related model, where the Verblunsky coefficients are given as follows. For $h\in C^{\omega}(\mathbb{T}^{d},\mathbb{R})$, $\lambda_1, \lambda_2 \in [0,1)$, we set
\begin{equation}\label{eq.sampling}
f(x,j)=\left\{\begin{aligned} \lambda_{1}e^{ i h(x)}\qquad \text{for } j=0 \mod 2,\\\lambda_{2}e^{i  h(x)}\qquad \text{for } j=1 \mod 2,\end{aligned}\right.
\end{equation}
and
\begin{equation}\label{eq.way}
\alpha_{n}(x) = f(x + n\omega, n), \; x \in \T^d , \; n \in \mathbb{Z}.
\end{equation}
Without loss of generality, we always assume that $\lambda_{1}^{2}+\lambda_{2}^{2}\neq 0$. For the resulting almost periodic CMV matrices, we have the following analog of Theorem~\ref{MainThm}:

\begin{Theorem}\label{MainThm1}
Let $\omega\in DC(\kappa,\tau)$ and $\rho>0$, $\lambda_{1}\in (0,1)$, then there exist constants $\epsilon=\epsilon(\kappa,\tau,\rho)$ and $r=r(\lambda_{1}) \in (0,1)$ such that for every $\lambda_2$ with $|\lambda_2|< r(\lambda_1)$ and every $h$ from a dense $G_\delta$ subset of
$$
\{ h \in C^{\omega}_{\rho}(\mathbb{T}^{d},\mathbb{R}) : \|h\|_{\rho}\leq \epsilon\},
$$
all gaps of the spectrum of the extended CMV matrix associated with the Verblunsky coefficients given by \eqref{eq.sampling}--\eqref{eq.way} are open, which in turn implies that the spectrum is a Cantor set in these cases.
\end{Theorem}

Of course when we say that all gaps are open, we are implicitly referring to the gap labelling associated with the base dynamics here. We leave the discussion of the relevant details to Section~\ref{AlmostModel}. The motivation for us to consider this particular model stems from two sources.

On the one hand, we are interested in studying more general almost periodic cases. A straightforward way to leave the setting of finite-dimensional tori (which describe the hulls occurring in quasi--periodic settings) is to take the product of a torus with a compact abelian group that is not also a finite-dimensional torus. The simplest such case is to take an additional factor of the form $\mathbb Z_p$. This case is naturally interesting because we can now model quasi-periodic perturbations of general periodic sequences, whereas before we were really just considering quasi-periodic perturbations of constant (or period-one) sequences. The general expectation is that the phenomena should be the same. However, the proofs become significantly more involved. For this very reason, we focus on the special case $p = 2$ here, as the resulting difficulties can still be dealt with via explicit calculations. We expect the general proof strategy to extend to larger values of $p$, but the reader will see from our proof that the calculations quickly become quite difficult.

On the other hand, the case $p = 2$ is of special interest because it is this extension of the quasi-periodic case that is needed to study quantum walks on the integer lattice with quasi-periodic coins. We describe this application below.

Consider the Hilbert space $\mathcal{H} = \ell^{2}(\mathbb{Z})\otimes\mathbb{C}^{2}$. A basis of $\mathcal{H}$ is given by elementary tensors $|n\rangle\otimes |\uparrow\rangle, |n\rangle\otimes |\downarrow\rangle$. Given the coins $$C_{n,t}=\begin{pmatrix}c_{n,t}^{11} &c_{n,t}^{12}\\c_{n,t}^{21}&c_{n,t}^{22}\end{pmatrix}\in \mathrm{U}(2)\qquad n,t\in \mathbb{Z},$$
the update rule of the quantum walk from time $t$ to time $t+1$ is as follows:
\begin{equation}\label{updateRule1}|n\rangle\otimes |\uparrow\rangle\mapsto c_{n,t}^{11}|n+1\rangle\otimes |\uparrow\rangle+c_{n,t}^{21}|n-1\rangle\otimes |\downarrow\rangle,\end{equation}
\begin{equation}\label{updateRule2}|n\rangle\otimes |\downarrow\rangle\mapsto c_{n,t}^{12}|n+1\rangle\otimes |\uparrow\rangle+c_{n,t}^{22}|n-1\rangle\otimes |\downarrow\rangle.\end{equation}
One can easily extend this rule to elements of $\mathcal{H}$ by linearity.

There are two cases of special interest:

$\bullet$ time-homogeneous: $C_{n,t}=C_{n}$ for every $t\in\mathbb{Z}$,

$\bullet$ space-homogeneous: $C_{n,t}=C_{t}$ for every $n\in\mathbb{Z}$.

We will focus our attention on the time-homogeneous case, that is, $C_{n}\in \mathrm{U}(2), n\in\mathbb{Z}$ are given. Then there is a unitary operator $U:\mathcal{H}\to\mathcal{H}$ such that the powers of $U$ provide the time evolution. Moreover, by choosing a basis of $\mathcal{H}$ in the following way:
$$\cdots,|-1\rangle\otimes|\uparrow\rangle,|-1\rangle\otimes|\downarrow\rangle,|0\rangle\otimes|\uparrow\rangle,|0\rangle\otimes|\downarrow\rangle,|1\rangle\otimes|\uparrow\rangle,|1\rangle\otimes|\downarrow\rangle\cdots,$$
the matrix representation of $U:\mathcal{H}\to\mathcal{H}$ is given by
\begin{equation}\label{eq.Unitary}U=\begin{pmatrix}
\cdots &\cdots&\cdots &\cdots&\cdots &\cdots&\cdots &\cdots\\
\cdots &0        &c_{-2}^{12}&0&0&0&0&\cdots\\
\cdots &c_{-1}^{21}&0&0&c_{-1}^{11}&0&0&\cdots\\
\cdots &c_{-1}^{22}&0&0&c_{-1}^{12}&0&0&\cdots\\
\cdots &0&0&c_{0}^{21}&0&0&c_{0}^{11}&\cdots\\
\cdots &0&0&c_{0}^{22}&0&0&c_{0}^{12}&\cdots\\
\cdots &0&0&0&0&c_{1}^{21}&0&\cdots\\
\cdots &\cdots&\cdots &\cdots&\cdots &\cdots&\cdots &\cdots
\end{pmatrix},\end{equation}
as can be readily checked using the update rules \eqref{updateRule1} and \eqref{updateRule2}. Recall that an extended CMV matrix with Verblunsky coefficients with odd index vanishing takes the form
$$\mathcal{E}=\begin{pmatrix}
\cdots &\cdots&\cdots &\cdots&\cdots &\cdots&\cdots &\cdots\\
\cdots &0        &-\alpha_{-4}&0&0&0&0&\cdots\\
\cdots &\overline{\alpha}_{-2}&0&0&\rho_{-2}&0&0&\cdots\\
\cdots &\rho_{-2}&0&0&-\alpha_{-2}&0&0&\cdots\\
\cdots &0&0&\overline{\alpha}_{0}&0&0&\rho_{0}&\cdots\\
\cdots &0&0&\rho_{0}&0&0&-\alpha_{0}&\cdots\\
\cdots &0&0&0&0&\overline{\alpha}_{2}&0&\cdots\\
\cdots &\cdots&\cdots &\cdots&\cdots &\cdots&\cdots &\cdots
\end{pmatrix},$$
where $\rho_{n}=\sqrt{1-|\alpha_{n}|^{2}}$. Then there is a unitary matrix
$$
\Lambda = \mathrm{diag} \, (\cdots,\lambda_{-1},\lambda_{0},\lambda_{1},\cdots)$$ such that $\Lambda^{*}U\Lambda=\mathcal{E},$
where $\mathcal{E}$ is the extended CMV matrix with associated Verblunsky coefficients given by
\begin{equation}\label{eq.VC}\alpha_{2n+1}=0,\quad\alpha_{2n+2}=\frac{\lambda_{2n+2}}{\lambda_{2n+1}}\overline{c}_{n+1}^{21},\quad n\in \mathbb{Z}.\end{equation}
The $\lambda_{n}$'s are specified as follows. With $U$ be given by \eqref{eq.Unitary}, write
$$
c_{n}^{kk}=|c_{n}^{kk}|e^{i\sigma_{n}^{k}}, n\in\mathbb{Z},k\in\{1,2\}, \sigma_{n}^{k}\in[0,2\pi)
$$
and define $\lambda_{n}$ by
$$
\begin{aligned}
&\lambda_{-1}=1\\
&\lambda_{0}=1\\
&\lambda_{2n+1}=e^{i\sigma_{n}^{2}}\lambda_{2n-1}\\
&\lambda_{2n+2}=e^{-i\sigma_{n}^{1}}\lambda_{2n}.
\end{aligned}
$$
One can consult \cite{CGMV} for this connection between quantum walks on the integer lattice and extended CMV matrices; it is commonly referred to as the CGMV connection.
Now as a consequence of Theorem \ref{MainThm1}, we have the following:

\begin{Theorem}\label{Thm.QW}
Consider the quantum walk with time-homogeneous coins given by $C_{n}\in \mathrm{SU}(1,1)$ and $c_{n}^{12}=\lambda e^{ih(x+n\omega)}$ with $\lambda\in (0,1), \omega\in DC(\kappa,\tau)$ and $h\in C_{\rho}^{\omega}(\mathbb{T}^{d},\mathbb{R})$ for some $\rho>0$. Then there exists $\epsilon=\epsilon(\kappa,\tau,\rho)$ such that for a generic $h$ in
$$\{h\in C_{\rho}^{\omega}(T^{d},\mathbb{R}):\|h\|_\rho<\epsilon\}$$
the spectrum of the corresponding unitary operator $U$ given by \eqref{eq.Unitary} has all gaps open. In particular, $U$ has Cantor spectrum in these cases.
\end{Theorem}

\subsection{Analyticity of Resonance Tongues}

The proofs of Theorem~\ref{MainThm} and Theorem \ref{MainThm1} may be facilitated by introducing an additional parameter, which plays the role of a coupling constant. Let us discuss this in the quasi-periodic case (i.e., Theorem~\ref{MainThm}), the discussion relevant to  Theorem \ref{MainThm1} is similar.

Basically, under explicit conditions (namely, the non-vanishing of Fourier coefficients) it can be shown that even if a gap is not open, it may be opened by wiggling this additional parameter. On the level of the associated Szeg\H{o} cocycles, we will then have two relevant parameters: the spectral parameter and the coupling constant.

To be specific, we will consider Verblunsky coefficients of the form
\begin{equation}\label{setting.var}
\alpha_{n}(x) = f(x+n\omega), \quad n\in \mathbb{Z},
\end{equation}
where
\begin{equation}\label{setting2.var}
f(x) = \lambda e^{i \delta h(x)},
\end{equation}
with $\omega, \lambda, h$ as above and $\delta \in \mathbb{R}$. Thus, relative to \eqref{setting}--\eqref{setting2}, we have just replaced $h$ by $\delta h$, and hence introduced a coupling constant in front of it. On the other hand, we will write the spectral parameter, which should be an element of the unit circle $\partial \mathbb{D}$, as $e^{i \theta}$ with $\theta \in \R$.

The Szeg\H{o} cocycle $(\omega, S(\theta,\delta))$ is then the quasi-periodic $\mathrm{SU}(1,1)$ cocycle over the base $x \mapsto x + \omega$ with the $(\theta,\delta)$-dependent map
$$
S(\theta,\delta) : \mathbb{T}^d \to \mathrm{SU}(1,1), \quad x \mapsto \frac{e^{-\frac{1}{2}i\theta}}{\sqrt{1-\lambda^{2}}}\begin{pmatrix} e^{i \theta} &-\lambda e^{-i\delta h(x)} \\ -\lambda e^{i\delta h(x)} e^{i \theta} & 1 \end{pmatrix}.
$$
The spectrum $\Sigma_\delta$ corresponds to the complement of the set of spectral parameters $e^{i\theta}$ for which the cocycle $(\omega, S(\theta,\delta))$ is uniformly hyperbolic. Denote the fibered rotation number of $S(\theta,\delta)$ by $\rot(\theta,\delta)$. Then, by the Gap Labelling Theorem, for each gap of $\Sigma_\delta$, there exists a $k\in\mathbb{Z}^{d}$ such that $\rot(\theta,\delta)=\frac{\langle k,\omega\rangle}{2} \mod \mathbb{Z}$ for all spectral parameters $e^{i\theta}$ in this gap. We say that this gap has \textit{label} $k$. One usually refers to the presence of this gap by saying that the ``gap with label $k$ is open.'' By contrast, if there is no actual gap of $\Sigma_\delta$ with label $k$, then one usually says that the ``gap with label $k$ is closed.'' The closed gap then refers to the unique point $e^{i \theta_k}$ in $\Sigma_\delta$ for which $\rot(\theta_k,\delta)=\frac{\langle k,\omega\rangle}{2} \mod \mathbb{Z}$ and one interprets that as the gap with label $k$ having degenerated to a single point, and hence closed.

Varying both parameters now and considering the preimage, we can define the set
\begin{equation}\label{eq.ResonanceTongue}
\mathcal{R}(k) = \left\{ (\theta,\delta)\in\mathbb{R}^{2} : \rot(\theta,\delta)=\frac{\langle k,\omega\rangle}{2} \mod \mathbb{Z} \right\},
\end{equation}
which is commonly called a \textit{resonance tongue}.

Note that for $\delta = 0$, the Verblunsky coefficients defined by \eqref{setting.var}--\eqref{setting2.var} take the constant value $\alpha_n(x) = \lambda$. The spectrum $\Sigma_{\delta = 0}$ is well known in this case and takes the form
$$
\Sigma_{\delta = 0} = \{e^{i\theta}:2\arcsin |\lambda| \leq \theta \leq 2\pi - 2 \arcsin |\lambda|\}
$$

Thus, only a single gap is open, namely the one corresponding to $k = 0$, while all others are closed. In particular, for each $k \not= 0$, the resonance tongue $\mathcal{R}(k)$ will have a \textit{tip} at $(\theta_{k,0},0)$, where $\theta_{k,0}$ is the unique choice for which $\rot(\theta_{k,0},0) = \frac{\langle k,\omega\rangle}{2} \mod \mathbb{Z}$.

A crucial step in our analysis is to show that under suitable assumptions, the boundaries of each $\mathcal{R}(k)$ are analytic, and that $\mathcal{R}(k)$ is non-degenerate (and in particular opens at the tip $(\theta_{k,0},0)$) if and only if the $k$-th Fourier coefficient of $h$ does not vanish. It then follows from $\hat h_k \not= 0$ that even if the $k$-th gap of the spectrum is closed for some value of the coupling constant, it opens as soon as the coupling constant is varied. In particular, requiring $\hat h_k \not= 0$ for all $k$ determines a dense $G_\delta$ subset of $h$'s for which the following is true: for every $k$ and all but countably many $\delta$, all gaps are open.

This discussion shows that a version of Theorem~\ref{MainThm} that is more reflective of how the result is obtained is given by the following theorem, and Theorem~\ref{MainThm} is an immediate consequence of it.

\begin{Theorem}\label{MainThm.refined}
Let $\omega\in \mathrm{DC}(\kappa,\tau)$, $\lambda \in (0,1)$, and $ \rho>0$ be fixed. Then there exists a constant $\epsilon=\epsilon(\kappa,\tau,\rho) > 0$ such that for every $h$ from a dense $G_\delta$ subset of
$$
\{ h \in C^{\omega}_{\rho}(\mathbb{T}^{d},\mathbb{R}) : \|h\|_{\rho}\leq \epsilon\},
$$
and all but countably many $\delta \in [-1,1]$, all gaps of the spectrum associated with the Verblunsky coefficients given by \eqref{setting.var}--\eqref{setting2.var} are open, which in turn implies that the spectrum is a Cantor set in these cases.
\end{Theorem}

The remainder of the paper is structured as follows. We begin with some preliminaries in Section~\ref{sec.2}, where we recall some standard notions and objects such as cocycles, the Lyapunov exponent, the rotation number, and reducibility. We provide an abstract criterion for analytic tongue boundaries in Section~\ref{section4}. With this criterion in hand we can then prove Theorem~\ref{MainThm.refined} in Section~\ref{sec.4} and Theorem~\ref{MainThm1} in Section~\ref{AlmostModel}. As explained above, Theorem~\ref{MainThm} follows from Theorem~\ref{MainThm.refined} and Theorem~\ref{Thm.QW} follows from Theorem~\ref{MainThm1}. We conclude the paper with two appendices that contain proofs of technical results that are needed in the main body of the paper.

\section*{Acknowledgment}
This paper is inspired by and borrows ideas from the paper \cite{Puig06} by Puig and Sim\'o. We are grateful to Joaquim Puig for some useful conversations about \cite{Puig06} and the philosophy underlying it.

\section{Preliminaries}\label{sec.2}

In this section we recall some fundamental concepts, which are important to our considerations below.

\subsection{Cocycles and the Lyapunov Exponent} \label{co-ly}
Suppose $X$ is a compact metric space, $T : X \to X$ is a homeomorphism, and $\nu$ is an ergodic Borel probability measure. Given $A \in C(X, \mathrm{SL}(2,\R))$, we associate the \textit{cocycle} $(T, A)$, which is given by the linear skew product
\begin{align*}
(T,A) :  X \times \R^2 & \to X \times \R^2 \\
(x,\phi) & \mapsto ( T x,A(x) \cdot \phi).
\end{align*}
For $n\in\mathbb{Z}$, $A_n$ is defined by $(T,A)^n=(T^n,A_n)$. Thus $A_{0}(x)=id$,
\begin{equation*}
A_{n}(x)=\prod_{j=n-1}^{0}A(T^{j}x)=A(T^{n-1}x)\cdots A(Tx)A(x) \text{ for } n \ge 1,
\end{equation*}
and $A_{-n}(x)=A_{n}(T^{-n}x)^{-1}$. The \textit{Lyapunov exponent} is defined by
$$
L(T,A) = \lim_{n \rightarrow \infty} \frac{1}{n} \int_{X} \ln \|A_{n}(x)\| \, d\nu(x).
$$
Here, $\| \cdot \|$ can of course be any norm on $\mathrm{SL}(2,\R)$, but if we choose the operator norm for convenience, then the existence of the limit follows immediately due to subadditivity.

In this paper, we will consider the following two special cases.
\begin{itemize}

\item  $X=\mathbb{T}^d$ and $T=R_{\omega}$ is the shift by $\omega \in X$, that is, $R_{\omega} x = x + \omega$. In this case $(\omega,A):=(R_{\omega},A)$ is called a \textit{quasi-periodic cocycle}. We generally assume that $\omega$ is such that $T$ is minimal, and in this case $\nu$ has to be normalized Lebesgue measure.

\item  $X=\mathbb{T}^d\times\mathbb{Z}_{\kappa}$, where $\kappa\in \Z^+$ and $T=T_{\omega}$, where $T_{\omega} (x,n)= (x+\omega,n+1)$. We can view $(T_{\omega},A)$ as an \textit{almost-periodic cocycle}.

\end{itemize}

%
%
%

\subsection{The Rotation Number}

Denote by $\mathbb S^{1}$ the set of unit vectors in $\mathbb{R}^{2}$ and consider a projective cocycle $F_{A}$ on $X\times\mathbb{S}^1$ defined by
$$
(x,\phi)\mapsto \Big( Tx, \frac{A(x)\phi}{\|A(x)\phi\|} \Big).
$$
If  $A\in C(X, \SL(2,\R))$ is homotopic to the identity, then there exists a lift $\tilde{F}_{A}$ of $F_{A}$ to $X \times \mathbb{R}$ such that $\tilde{F}_{A}(x,\phi) = (Tx, \tilde{f}_{A}(x,\phi))$,  where $\tilde{f}_{A} :X \times \mathbb{R}\rightarrow\mathbb{R}$ is a continuous lift with the following properties:
\begin{itemize}
  \item $\tilde{f}_{A}(x,\phi + 1) = \tilde{f}_{A}(x,\phi) + 1;$
  \item for every $x \in X, \tilde{f}_{A}(x, \cdot) : \mathbb{R} \rightarrow \mathbb{R}$ is a strictly increasing homeomorphism;
  \item if $\pi_2$ is the projection map  $X \times \mathbb{R}\rightarrow X \times \mathbb{S}^1:(x,\phi) \mapsto (x, e^{2\pi i\phi})$, then $F_{A} \circ \pi_2 = \pi_2 \circ \tilde{F}_{A}$.
\end{itemize}
If $(X, T)$ is uniquely ergodic, then the number
\begin{equation*}
\rho(T,A) = \lim_{n\rightarrow \infty}\frac{\tilde{f}_{A}^{n}(x,\phi) -\phi}{n} \mod \Z
\end{equation*}
is independent of  $(x,\phi) \in X \times \mathbb{R/Z}$ and the choice of the lift of $F_{A}$. It is called the {\it fibered rotation number} of $(T,A)$; see \cite{H,JM82} for details.


Once we have the definition of the fibered rotation number for $\mathrm{SL}(2,\mathbb{R})$ cocycles, we can parlay this concept to $\mathrm{SU}(1,1)$ cocycles, since
$$
M=\frac{1}{1+i}\begin{pmatrix}1 &-i\\ 1 &i\end{pmatrix}
$$
induces an isomorphism between $\mathrm{SL}(1,1)$ and $\mathrm{SL}(2,\mathbb{R})$:
$$
M^{-1} \mathrm{SU}(1,1) M = \mathrm{SL}(2,\mathbb{R}).
$$

\subsection{Reducibility}

Recall that in the quasi-periodic setting we have $X=\mathbb{T}^d$, $T = R_{\omega}$, and $A_1 \in C^\omega(\mathbb{T}^d, \mathrm{SL}(2,\R))$.

The cocycle $(\omega,A_1)$ is said to be \textit{analytically conjugate} to $(\omega, A_2)$ if there exists $B \in C^{\omega}(2\T^{d}, \mathrm{SL}(2,\R))$ such that
$$
B^{-1}(\cdot+\omega) A_1(\cdot) B(\cdot) = A_2(\cdot).
$$
The cocycle $(\omega, A)$ is said to be \textit{almost reducible} if the closure of its analytic conjugations contains a constant.

We also need the following well-known result of Eliasson \cite{Eli92}:

\begin{Theorem} \cite{Eli92}\label{Thm.reduci}
Let $r > 0$, $\omega \in \mathrm{DC}(\kappa,\tau)$, and $A_{0}\in \mathrm{SL}(2,\mathbb{R})$. Then there is a constant $\epsilon_{0} = \epsilon_{0}(\kappa,\tau,r,\| A_{0}\|)$ such that if $A \in C^{\omega}_{r}(\mathbb{T}^{d},\mathrm{SL}(2,\mathbb{R}))$ is real analytic with
$$
\| A - A_{0} \|_{r} \leq \epsilon_{0}
$$
and the rotation number of the cocycle $(\omega, A)$ satisfies
$$
\|2\rho(\omega,A) - \langle n, \omega \rangle\|_{\mathbb{R}/\mathbb{Z}} \geq \frac{\kappa}{|n|^{\tau}}, \quad \forall  0\neq n\in \mathbb{Z}^{d}
$$
 or $2\rho(\omega,A) = \langle n, \omega \rangle\mod \Z$ for some $n \in \mathbb{Z}^{d}$, then $(\omega,A)$ is analytically reducible to a constant, that is, there exists $Z(x)\in C^{\omega}(2\mathbb{T}^{d},\mathrm{SL}(2,\mathbb{R}))$ such that
$$Z^{-1}(x+\omega)A(x)Z(x)=B,$$
where $B\in\mathrm{SL}(2,\mathbb{R})$ is a constant matrix.
\end{Theorem}
\begin{Remark}\label{uniform}
If $A_{0}$ varies in a bounded set, we can choose  $\epsilon_0= \epsilon_0(\kappa,\tau,r) $ to be independent of $\|A_{0}\|$.
\end{Remark}

\section{A Criterion for Analytic Tongue Boundaries}\label{section4}

In this section we give a general criterion for a parametrized family of cocycles to have analytic tongue boundaries.

We first describe the basic setting. Consider analytic quasi-periodic  cocycles $(\omega, A_{\Gamma}(\cdot))$ depending on the parameters $\omega \in \mathrm{DC}(\kappa,\tau)$ and $\Gamma=(\theta,\delta)\in \mathbb{R}^{2}$. In typical applications we will take $\theta$ to be the spectral parameter and $\delta$ to be the coupling constant.

Let $\Gamma_{0}=(\theta_{0},\delta_{0})$ be fixed and, with a neighborhood $N(\Gamma_{0})$ of $\Gamma_{0}$, assume the following:
\begin{itemize}
\item[(H1):] $2 \rho(\omega,A_{\Gamma_{0}}(\cdot))=\langle k,\omega\rangle\mod\mathbb{Z}$ and is a non-decreasing (or non-increasing) function of $\theta$ in a neighborhood of $\theta_{0}$.
\item[(H2):]  $(\omega,A_{\Gamma_{0}}(\cdot))$ is not uniformly hyperbolic.
\item[(H3):]  $\sup_{\Gamma\in N(\Gamma_{0})} \|A_{\Gamma}(\cdot)-A_{\theta}\|_{r}\leq\epsilon_{0}(\kappa,\tau,r)$ for some constant $A_{\theta}\in\SL(2,\mathbb{R})$,
\end{itemize}
where $\epsilon_0= \epsilon_0(\kappa,\tau,r)$ is defined in Theorem \ref{Thm.reduci} (see also Remark \ref{uniform}).
By the assumptions (H1) and (H3), along with $\omega \in \mathrm{DC}(\kappa,\tau)$, Theorem \ref{Thm.reduci} implies that $(\omega,A_{\Gamma_{0}}(x))$ is analytically reducible, that is, there exists some $Z(x)\in C^{\omega}(2\mathbb{T}^{d},\SL(2,\mathbb{R}))$ such that
\begin{equation}\label{red}Z^{-1}(x+\omega)A_{\Gamma_{0}}(x)Z(x)=B.\end{equation}
By (H2), the constant  matrix $B$ is in fact parabolic, that is, we can rewrite \eqref{red} as
$$
Z^{-1}(x+\omega)A_{\Gamma_{0}}(x)Z(x)=\begin{pmatrix}1 &c\\0 &1\end{pmatrix}.
$$
Then we have for $\Gamma$ in a neighborhood of $\Gamma_{0}$,
\begin{equation}\label{eq.newconj}
Z^{-1}(x+\omega)A_{\Gamma}(x)Z(x)=\begin{pmatrix}1 &c\\0 &1\end{pmatrix}e^{P(x,\zeta)+O_{2}(x,\zeta)}
\end{equation}
where $\zeta = \Gamma-\Gamma_{0}$, $P(x,\zeta)$ is a linear term in $\zeta$, and $O_{2}(x,\zeta)$ is the collection of higher order terms in $\zeta$.  Let
$$
[P]=\begin{pmatrix} [P_{11}]& [P_{12}]\\ [P_{21}]&- [P_{11}]\end{pmatrix}
$$
be the average of $P(x,\zeta)$ with respect to $x$ over $\T^{d}$, that is, the integral with respect to normalized Lebesgue measure. We have the following result:

\begin{Proposition}\label{twoCases}
Let $\omega \in \mathrm{DC}(\kappa,\tau), \Gamma_{0}=(\theta_{0},\delta_{0})$ be fixed and assume that for a suitable neighborhood of $\Gamma_{0}$, the hypotheses {\rm (H1), (H2)} and {\rm (H3)} are satisfied. Then we have the following:
\begin{enumerate}
\item If $c\neq 0$ and  the coefficient of $(\theta-\theta_0)$ in $[P_{21}]$ is strictly non-zero,
then the tongue boundary $\theta=\theta(\delta)$ is real analytic in a neighborhood of $\delta_{0}$.
\item  If $c=0$ and  the coefficient of $(\theta-\theta_0)^{2}$ in the expression of $\det [P]$ is nonzero, then the two tongue boundaries $\theta_{i}=\theta_{i}(\delta),i=1,2$ with $\theta_{i}(\delta_{0})=\theta_{0}$ are real analytic in a neighborhood of $\delta_{0}$.
\end{enumerate}
\end{Proposition}

\begin{Remark}
This criterion was first proved for the quasi-periodic Hill equation (a special linear quasi-periodic differential equation) by Puig-Sim\'o \cite{Puig06}. We generalize it to {\it general} $\SL(2,\R)$-valued analytic quasi-periodic cocycles.
\end{Remark}

Before presenting the proof of Proposition~\ref{twoCases}, we first need to introduce a useful tool in the following subsection.

\subsection{Existence of Counterterms}

In this subsection we give two criteria for the existence of counterterms, which make the original near constant cocycle reducible to the same constant cocycle.

We first introduce the following definition, which is inspired by \cite{Puig06}:

\begin{Definition}\label{admi}
Let $A_{0}\in \sl(2,\mathbb{R})$ and $C,S:  \sl(2,\mathbb{R}) \to \sl(2,\mathbb{R}) $ with $C^{2}=C$. We say the quartet $(A_{0},C,S,\omega)$ is admissible if there exist positive constants $c',\nu$ such that for all analytic $P\in C_{\rho}^{\omega}(\mathbb{T}^{d}, \sl(2,\mathbb{R}))$, the equations
\begin{equation}\label{eq.homoeq}\begin{aligned}
 -e^{-A_{0}}X(x+\omega)e^{A_{0}}+X(x)=-(P-C([P])),\\
[X]=S([P])\end{aligned}\end{equation}
have a unique real analytic solution $X:2\mathbb{T}^{d}\to \sl(2,\mathbb{R})$ that satisfies the estimate
$$\|X\|_{\rho-r}\leq c'\frac{\|P\|_{\rho}}{r^{\nu}}$$ for any $0<r<\rho$,
where $[\cdot]$ denotes the average with respect to normalized Lebesgue measure.
\end{Definition}

In the following, let $\mathcal{O}$ be a set of external multi-parameters; an element of $\mathcal{O}$ will be denoted by $\zeta$. We have the following theorem.

\begin{Theorem}\label{KAM}
Let $A_{0}=\begin{pmatrix}0 & c \\0&0\end{pmatrix}$ and assume that $(A_{0},C,S,\omega)$ is admissible with positive constants $c',\nu$. Let $\rho_{0}>0$ and $p\in \mathbb{Z}_{+}$, then there exists $\epsilon = \epsilon(A_{0},\rho_{0},c',\nu) > 0$ such that for any $f\in C^{\omega}_{\rho_{0}}(\mathbb{T}^{d},\sl(2,\mathbb{R}))$ with $\|f\|_{\rho_{0}}\leq\epsilon$ and any $|\chi|\leq1$, there exists $\xi^{*}\in\sl(2,\mathbb{R})$ with $C(\xi^{*})=\xi^{*}$ such that the cocycle  $$(\omega, e^{\chi^{p}A_{0}}e^{\chi^{p+1} f(x)}e^{-\chi^{p+1}\xi^{*}})$$ is analytically conjugated into $(\omega,e^{\chi^{p}A_{0}})$ by a conjugation $Z(x)=e^{\chi X(x)}$, where $X\in C^{\omega}(2\mathbb{T}^{d},\sl(2,\mathbb{R}))$ with $$\|X\|_{\rho_{0}/2}\leq \tilde{c}(\nu)\epsilon\rho_{0}^{-\nu},$$ where $\tilde{c}(\nu)$ is a constant. Moreover, if $P$ and $\chi$ depend analytically on $\zeta\in\mathbb{R}^{2}$ in a neighborhood around the origin, then $X$ and $\xi^{*}$ depend analytically on $\zeta$ in a smaller neighborhood.
\end{Theorem}

\begin{Remark}
The idea of the existence of a counterterm was first introduced by Moser \cite{Moser}; one can consult \cite{BMS,Katok,Kri,Puig06} for further developments.
\end{Remark}

We also introduce another version of a counterterm existence result, which is slightly different due to the non-commutativity property of $\SL(2,\R)$:

\begin{Theorem}\label{KAM1}
Let $A_{0}=\begin{pmatrix}0 &c\\0&0\end{pmatrix}$ and assume that $(A_{0},C,S,\omega)$ is admissible with positive constants $c',\nu$. Let $\rho_{0}>0$ and $p\in\mathbb{Z}_{+}$, then there exists $\epsilon = \epsilon(A_{0},\rho_{0},c',\nu) > 0$ such that for any $f\in C^{\omega}_{\rho_{0}}(\mathbb{T}^{d},\sl(2,\mathbb{R}))$ with $\|f\|_{\rho_{0}}\leq\epsilon$ and any $|\chi|\leq1$, there exists $\xi^{*}\in\sl(2,\mathbb{R})$ with $C(\xi^{*})=\xi^{*}$ such that the cocycle
$$
(\omega, e^{\chi^{p}(A_{0}+\chi f(x)-\chi\xi^{*})})
$$
is analytically reducible to $(\omega, e^{\chi^{p}A_{0}})$ by a conjugation $Z(x)=e^{\chi X(x)}$, where $X\in C^{\omega}(2\mathbb{T}^{d},\sl(2,\mathbb{R}))$ with $$\|X\|_{\rho_{0}/2}\leq \tilde{c}(\nu)\epsilon\rho_{0}^{-\nu},$$ where $\tilde{c}(\nu)$ is a constant. Moreover, if $P$ and $\chi$ depend analytically on $\zeta\in\mathbb{R}^{2}$ in a neighborhood around the origin, then $X$ and $\xi^{*}$ depend analytically on $\zeta$ in a smaller neighborhood.
\end{Theorem}

The proof of Theorem \ref{KAM1} is almost the same as that of Theorem \ref{KAM} and we present these proofs in Appendix~\ref{appendix1}.

\subsection{Proof of Proposition \ref{twoCases}:}

Let $\zeta=\Gamma-\Gamma_{0}$, and denote $\beta=\delta-\delta_{0}$,  $\eta=\theta-\theta_0$, and $\zeta=(\eta,\beta)$.  It is readily checked that if  $A_{0}=\begin{pmatrix}0 & c \\0&0\end{pmatrix}$ and $\omega \in \mathrm{DC}(\kappa,\tau)$,  then there exist choices of $C$ and $S$ such that the quartet $(A_{0},C,S,\omega)$  is admissible. Thus for $\Gamma$ in a suitable small neighborhood of $\Gamma_{0}$, one can apply Theorem \ref{KAM},  and  there exists a counterterm $\xi^{*}(\zeta)\in \sl(2,\mathbb{R})$ such that the cocycle
\begin{equation}\label{eq.newsys} \left(\omega, \begin{pmatrix} 1 &c\\ 0 &1\end{pmatrix}e^{P(x,\zeta)+O_{2}(x,\zeta)}e^{-\xi^{*}(\zeta)}\right)
\end{equation}
is analytically reducible to $\left(\omega, \begin{pmatrix}1 &c\\ 0 &1\end{pmatrix}\right)$.
We prove the proposition by distinguishing two different cases:

\smallskip
\textbf{Case 1 $c\neq 0$:}  In this case, let $A_{0}=\begin{pmatrix}0 &c\\ 0 &0\end{pmatrix}$ and take $\begin{pmatrix}1 &0\\ 0 &-1\end{pmatrix}, $ $\begin{pmatrix}0 &1\\ 0 &0\end{pmatrix},$ $\begin{pmatrix}0 &0\\1 &0\end{pmatrix}$ as a basis of  $\sl(2,\mathbb{R})$, let $[X]=\begin{pmatrix}X_{11}&X_{12}\\X_{21}&-X_{11}\end{pmatrix}$. Define the linear operator $ad_{A_{0}}: \sl(2,\mathbb{R})\to\sl(2,\mathbb{R}) $,
$$
ad_{A_{0}}([X])=e^{-A_{0}}[X]e^{A_{0}}-[X] = \begin{pmatrix} -cX_{21} &2cX_{11}-c^{2}X_{21}\\ 0 &cX_{21} \end{pmatrix}.
$$
The spectrum of $ad_{A_{0}}$ reduces to the zero eigenvalue, its kernel is the linear subspace spanned by $\begin{pmatrix}0 &1\\0 &0\end{pmatrix}$. Note that the equation \eqref{eq.homoeq} can be solved by comparing Fourier coefficients if and only if $ad_{A_{0}}\circ S(\cdot)+C(\cdot)=Id$.
Therefore, we can choose
$$
C([X]) = \begin{pmatrix} 0 & 0 \\ X_{21} & 0 \end{pmatrix}, \; S([X]) = \begin{pmatrix} -\frac{cX_{11} - X_{12}}{2c} & 0 \\ -\frac{X_{11}}{c} & \frac{cX_{11} - X_{12}}{2c} \end{pmatrix}
$$
to make the quartet $(A_{0},C,S,\omega)$ admissible.

Since the counterterm given by Theorem \ref{KAM} satisfies $C(\xi^{*})=\xi^{*}$, we have
$$
\xi^{*}(\zeta) = \begin{pmatrix}0 &0\\  \xi^{*}_{21}(\zeta) & 0 \end{pmatrix},
$$
and for the values of $\zeta$ for which $ \xi^{*}_{21}(\zeta)=0$, the  cocycles $(\omega, A_{\Gamma})$ with parameters $(\theta(\delta),\delta)$ are reducible to $\left(\omega, \begin{pmatrix}1 &c\\ 0 &1\end{pmatrix}\right)$. Hence to prove item $(1)$,  we only need to show that the equation
$$
\xi^{*}_{21}(\zeta) = \xi^{*}_{21}(\eta,\beta) = 0
$$
can be inverted to obtain an analytic function $\eta=\eta(\beta)$. According to \eqref{eta} and the notations there, we can set $\eta_{\infty}=0$ and obtain $\xi^{*}=\eta_{0}$ recursively, in particular, as $c\neq 0$, we have
$$
\xi^{*}_{21}(\eta,\beta)=[P_{21}(\eta,\beta)]+O_{2}(\eta,\beta).
$$
Thus under condition (1) of Proposition~\ref{twoCases}, the coefficient of $\eta$ in $[P_{21}]$ is nonzero, and hence the Implicit Function Theorem gives rise to a real analytic function $\eta=\eta(\beta)$. This finishes the proof of the case $c\neq 0$.

\smallskip
\textbf{Case 2 $c=0$:} In this case $A_{0}=0$, and taking $C=Id$ and $S=0$ makes the quartet $(A_{0},C,S,\omega)$ admissible. Therefore, the condition $C(\xi^{*})=\xi^{*}$ does not provide extra information, so we need to take further  averaging steps. The idea is to transform the question into Case 1.

We first point out that after $r$  averaging   steps,    we obtain the expansion of $\xi^{*}(\zeta)$ up to order $r$ in the following form:
\begin{equation}\label{eq.newexpansion}
\begin{pmatrix}S_{3}^{(r)} &S_{2}^{(r)}\\-S_{1}^{(r)} &-S_{3}^{(r)}\end{pmatrix}+O_{r+1}(\zeta),
\end{equation}
where $S_{i}^{(r)}(\zeta)=\sum_{1\leq s\leq r}D_{s,i}^{(r)}(\zeta)$ for $i=1, 2, 3$ and the polynomials $D_{s,i}^{(r)}$ for $i=1,2,3$ have the following properties:

(i) $D_{s,i}^{(r)}(\zeta)$ are homogeneous of degree $s$ in $\zeta$;

(ii) If $s<t\leq r$, then $D_{s,i}^{(t)}=D_{s,i}^{(s)}$.\newline
Moreover, one can explicitly calculate an expression for $D_{1,i}^{1}$, $i=1,2,3$, by averaging once. Theorem~\ref{KAM} implies that the cocycle $$
\left(\omega, e^{P(x,\zeta)+O_{2}(x,\zeta)}e^{-\xi^{*}(\zeta)}\right)
$$
is analytically reducible to $(\omega, Id)$ via a conjugation
$B$ that is close to the identity. That is, there exists $B(x)\in C^{\omega}(2\mathbb{T}^{d},\mathrm{SL}(2,\mathbb{R}))$ such that
$$
B^{-1}(x+\omega)e^{P(x,\zeta)+O_{2}(x,\zeta)}e^{-\xi^{*}(\zeta)}B(x)=Id,
$$
therefore,
$$
B^{-1}(x+\omega)e^{P(x,\zeta)+O_{2}(x,\zeta)}B(x)=B^{-1}(x)e^{\xi^{*}(\zeta)}B(x).
$$
This implies that the original cocycles are now conjugated into a new cocycle $(\omega, B^{-1}(x)e^{\xi^{*}(\zeta)}B(x))$.

We need the following lemma:

\begin{Lemma}\label{nonzero}
If $\det \xi^{*}(\zeta)>0$, then the rotation number
$\rho(\omega, e^{P(x,\zeta)+O_{2}(x,\zeta)}) $
 is strictly non-zero.
\end{Lemma}

\begin{pf}
Since the conjugation $B$ is close to the identity, thus with zero degree, it is enough for us to prove that
$$
\rho(\omega, B^{-1}(x)e^{\xi^{*}(\zeta)}B(x)) \neq 0.
$$

Consider the following problem: let $u_{n},v_{n}$ be the sequences satisfying the following iteration:
$$
\begin{pmatrix}u_{n+1}\\ v_{n+1}\end{pmatrix}=\begin{pmatrix}a_{11} & a_{12}\\a_{21}&a_{22}\end{pmatrix}\begin{pmatrix}u_{n}\\v_{n}\end{pmatrix},
$$
where $A(\cdot)=\begin{pmatrix}a_{11}(\cdot) & a_{12} (\cdot)\\a_{21}(\cdot)&a_{22}(\cdot)\end{pmatrix}\in\SL(2,\mathbb{R})$,
and let $R_{n}e^{i\psi_{n}}=u_{n}+iv_{n}$ with $0\leq |\psi_{n+1}-\psi_{n}|<\pi$. Then we have
$$
\begin{aligned}R_{n+1}e^{i\psi_{n+1}}&=u_{n+1}+iv_{n+1}\\
&=(a_{11}u_{n}+a_{12}v_{n})+i(a_{21}u_{n}+a_{22}v_{n})\\
&=R_{n}(a_{11}\cos\psi_{n}+a_{12}\sin\psi_{n})+iR_{n}(a_{21}\cos\psi_{n}+a_{22}\sin\psi_{n}).
\end{aligned}
$$
Therefore we have
$$
\begin{aligned}
\frac{R_{n+1}}{R_{n}}e^{i(\psi_{n+1}-\psi_{n)}}=\frac{(a_{11}\cos\psi_{n}+a_{12}\sin\psi_{n})+i(a_{21}\cos\psi_{n}+a_{22}\sin\psi_{n})}{e^{i\psi_{n}}},
\end{aligned}
$$
and taking the imaginary part of both sides, we obtain
$$
\frac{R_{n+1}}{R_{n}}\sin(\psi_{n+1}-\psi_{n})=a_{21}\cos^{2}\psi_{n}+(a_{22}-a_{11})\sin\psi_{n}\cos\psi_{n}-a_{12}\sin^{2}\psi_{n}.
$$

Denote $K=\begin{pmatrix}a_{21} &\frac{a_{22}-a_{11}}{2}\\ \frac{a_{22}-a_{11}}{2}&-a_{12}\end{pmatrix}$. If $K$ is a definite quadratic form, then $\psi_{n+1}-\psi_{n}$ has fixed sign, which in turn implies that its rotation number $\rho(\omega,A)$ is strictly non-zero. Note that
$$
\det K=-a_{12}a_{21}-\frac{(a_{11}-a_{22})^{2}}{4}=\frac{4-(a_{11}+a_{22})^{2}}{4}.
$$
Thus $K$ is definite if and only if the trace of $A$ satisfies $|\tr(A)|<2$. Let $A=B^{-1}(x)e^{\xi^{*}(\zeta)}B(x)$. Then $\det \xi^{*}(\zeta) > 0$ implies that $|\tr(A)|<2$, and the lemma follows.
\end{pf}

According to \eqref{eq.newexpansion} and  \eqref{eta} with $c=0$, we have $\xi^{*}=[P]+O_{2}(\zeta)$, under the condition (ii), we can apply the Weierstrass Preparation Theorem to $\det\xi^{*}$ and obtain the following:
$$\det\xi^{*}(\eta,\beta)=F(\eta,\beta)(\eta^{2}+g_{1}(\beta)\eta+g_{2}(\beta)),$$
where $g_{1},g_{2}$ and $F$ are real analytic functions with $F(0,0)\neq0$.
Fix $\beta=0$, by Lemma \ref{nonzero} and (H1), the local non-decreasing (or non-increasing) property of the rotation number with respect to $\eta$, the rotation number changes sign when $\eta$ crosses zero. Then, according to the continuity of the rotation number, for each nonzero $\beta$ in a neighborhood of zero, we can find a suitable real $\eta^{*}=\eta^{*}(\beta)$ such that the rotation number is zero. This yields the following expression,
$$
\det\xi^{*}(\eta,\beta)=F(\eta,\beta)(\eta-\eta^{*}_{1}(\beta))(\eta-\eta^{*}_{2}(\beta)),
$$
where $\eta^{*}_{1},\eta^{*}_{2}$ are real analytic functions (see \cite{Rel69}). Since $\eta_{1}^{*},\eta_{2}^{*}$ both vanish at zero,
there are two possible subcases:

\smallskip
\textbf{Case 2-1:} If $\eta^{*}_{1}(\beta),\eta^{*}_{2}(\beta)$ coincide in a neighborhood of $0$, then their analyticity implies $\eta^{*}_{1}(\beta)=\eta^{*}_{2}(\beta)$ for all $\beta$, and they both give the tongue boundary.

\smallskip

\textbf{Case 2-2:} There exist a non-zero constant $C$ and a positive integer $p$ such that $$\eta^{*}_{1}(\beta)-\eta^{*}_{2}(\beta)=C\beta^{p}+O_{p+1}(\beta).$$
 Since we are interested in the roots of $\det\xi^{*}(\zeta)=0$, we can scale $\xi^{*}$ such that $\xi_{12}^{*}(\zeta),\xi_{21}^{*}(\zeta)$ start with $\eta+\cdots$.
Moreover, under condition (ii), we may assume $\eqref{eq.newexpansion}$ to take the following form,
\begin{equation}\label{eq.normalform}\left\{\begin{aligned}
S^{(r)}_{2}=\eta+\sigma_{1}(\beta)+\eta\rho_{1}(\eta,\beta)\\
-S^{(r)}_{1}=\eta+\sigma_{2}(\beta)+\eta\rho_{2}(\eta,\beta)\\
S^{(r)}_{3}=\qquad \sigma_{3}(\beta)+\eta\rho_{3}(\eta,\beta)
\end{aligned}\right. ,
\end{equation}
where $\sigma_{i}(\beta)=\sum_{k\geq p}m_{i,k}\beta^{k},i=1,2,3, (m_{1,p}-m_{2,p})^{2}+m_{3,p}^{2}>0$ and $\eta$ is a suitable redefinition which contains terms of $\beta$ with order less than $p$, $\rho_{j}$ are polynomials with maximal degree $r-1$. This conclusion needs some explanation:

Recall \eqref{eq.newexpansion}, since the coefficient of term of degree 1 in $\eta$ in $S^{(r)}_{1}$ equals 1, there is a conjugation induced by $\begin{pmatrix}1 &-c_{3} \\0 &1\end{pmatrix}$, where $c_{3}$ denotes the coefficient of $\eta$ in $S_{3}^{(r)}(\zeta)$, making $S^{(r)}_{3}$ contain no linear $\eta$ term. Therefore, we may assume that
$$
\left\{\begin{aligned}
S_{2}^{(r)}=\eta+\sigma_{1}(\beta)+\eta\rho_{1}(\eta,\beta)\\
-S_{1}^{(r)}=\eta+\sigma_{2}(\beta)+\eta\rho_{2}(\eta,\beta)\\
S_{3}^{(r)}=\qquad \sigma_{3}(\beta)+\eta\rho_{3}(\eta,\beta)
\end{aligned}\right. ,
$$
where $\eta$ is a suitable redefinition and contains the terms of order less than $p$ in $\beta$ and $\sigma_{j}$ are the polynomials in $\beta$ of maximal degree $r$.
Let $k_{j}$ be the minimal degree of $\sigma_{j}$ if $\sigma_{j}\neq 0$, otherwise, take $k_{j}=\infty$, define $k=\min \{k_{1},k_{2},k_{3}\}\leq r$. Imposing a change of variables in $\eta$ by $\eta=\gamma\beta^{k}$, then $S^{(r)}_{1}S_{2}^{(r)}-(S_{3}^{(r)})^{2}$ becomes
$$
\beta^{2k}(\gamma^{2}+(m_{1,k}+m_{2,k})\gamma+m_{1,k}m_{2,k}-m_{3,k}^{2}+O(\beta)),
$$
where $m_{j,k}$ is the coefficient of $\beta^{k}$ in $\sigma_{j}$ for $j=1,2,3$. By factoring out $\beta^{2k}$ and neglecting $O(\beta)$ terms, we obtain the following equation,
\begin{equation}\label{eq01}
\gamma^{2}+(m_{1,k}+m_{2,k})\gamma+m_{1,k}m_{2,k}-m_{3,k}^{2}=0
\end{equation}
with discriminant
$$
(m_{1,k}-m_{2,k})^{2}+4m_{3,k}^{2}.
$$
Therefore, \eqref{eq01} has one multiple root if and only if $m_{1,k}=m_{2,k},m_{3,k}=0$. In this case, we put the term $m_{i,k}\beta^{k}$ into $\eta$ of the corresponding entry of $\xi^{*}(\zeta)$, therefore, it becomes $\hat{\eta}$ and we still use the same notation $\eta$, and the minimal degree of $\sigma_{j}$ increases. This process can be done until $k=p$, since in this case, $$\eta^{*}_{1}(\beta)-\eta^{*}_{2}(\beta)=C\beta^{p}+O_{p+1}(\beta)$$ with $C\neq 0$, which means \eqref{eq01} has two different roots for $k=p$. On the other hand, the process will not stop for any $k<p$, if so, we have $$\eta^{*}_{1}(\beta)-\eta^{*}_{2}(\beta)=\tilde{C}\beta^{k}+O_{k+1}(\beta),$$ which contradicts the assumption that $\eta^{*}_{1}(\beta)$ and $\eta^{*}_{2}(\beta)$ have contact of order $p$. Let $\hat{\eta}=\eta-P(\beta)$, where $P$ is a suitable polynomial of degree at most $p$ in $\beta$. With this redefinition of $\eta$, we arrive at \eqref{eq.normalform}.

With the form of \eqref{eq.normalform}, let $\gamma_{1}$ be one of the roots of the following equation,
\begin{equation}\label{eq.equation1}
\gamma^{2}-(m_{1,p}+m_{2,p})\gamma+m_{1,p}m_{2,p}-m_{3,p}^{2}=0,
\end{equation}
and impose the change of variables $\eta=\gamma_{1}\beta^{p}+\Delta\beta^{p}$. Then $B^{-1}(x)e^{\xi^{*}(\zeta)}B(x)$ becomes
\begin{equation}\label{eq.newnormalform}
\exp\left(\beta^{p}\left(\begin{pmatrix}m_{3,p} &(m_{1,p}+\gamma_{1})+\Delta\\ -(m_{2,p}+\gamma_{1})-\Delta &-m_{3,p}\end{pmatrix}+\beta R(x,\zeta)\right)\right).
\end{equation}
Here we used the fact that $X$ is of order $O(\beta)$.
After one step of conjugation induced by
$$
\begin{pmatrix}\frac{1}{\sqrt{b+\Delta}} &0\\ \frac{a}{\sqrt{b+\Delta}} &\sqrt{b+\Delta}\end{pmatrix},
$$
where $a=m_{3,p}, b=m_{1,p}+\gamma_{1},c=m_{2,p}+\gamma_{1}$ with $a^{2}=bc$, since $b,c$ cannot be zero at the same time (otherwise the order of contact of $\eta^{*}_{1}$ and $\eta^{*}_{2}$ will be $p+1$), we  may assume that $b>0$, the case $b<0$ can be dealt similarly, the cocycle becomes
$$
\exp\left(\beta^{p}\left(\begin{pmatrix}0 &1\\ -\Delta(b+c+\Delta) &0\end{pmatrix}+\beta^{2}R_{1}(x,\Delta,\beta)\right)\right).
$$
The extra $\beta$ factor is imposed to make sure the perturbation term is sufficiently small, and can be obtained by pushing the averaging one more step forward. We use Theorem~\ref{KAM1} for $\chi=\beta, A_{0}=\begin{pmatrix}0 &1\\0 &0\end{pmatrix}$ to obtain the counterterm $\xi^{*}(\Delta,\beta)$ such that $(\omega, e^{\beta^{p}(A_{0}+\beta(P(x,\Delta,\beta)-\xi^{*}(\Delta,\beta)))})$ is analytically conjugated into $(\omega, e^{\beta^{p}A_{0}})$. The admissible assumption associated with the case $c\neq0$ requires $\xi^{*}$ to take the form $\begin{pmatrix}0 &0\\ \xi^{*}_{21} &0\end{pmatrix}$, where $\xi^{*}_{21}(\Delta,\beta)$ is a real analytic function. Therefore, the equation
\begin{equation}\label{eq.tb}
-\beta\xi^{*}_{21}(\Delta,\beta)+\Delta(b+c+\Delta)=0
\end{equation}
describes a branch of the tongue boundary. Note also that under our assumption $b>0$, we have $b+c>0$ and therefore we can deduce from \eqref{eq.tb} that $\Delta=\Delta(\beta)=O(\beta)$ and therefore $$\eta_{1}(\beta)=\gamma_{1}\beta^{p}+\Delta(\beta)\beta^{p}=\gamma_{1}\beta^{p}+O(\beta^{p+1}).$$ The other branch, which corresponds to $\gamma_{2}$, the other root of equation \eqref{eq.equation1}, can be treated similarly. This completes the proof of the analyticity of the tongue boundary. \qed

\subsection{Transitivity at the Origin}

Proposition \ref{twoCases} gives two cases when the tongue boundaries are analytic. In order to show generic Cantor spectrum, one still needs to show that the two tongue boundaries of a certain resonance have a transversality at the origin (say $\delta_{0}=0$).  On the other hand, while the general criterion for analytic boundaries can be better presented in  $\SL(2,\R)$, we found that when it comes to concrete  verification, it facilitates the calculations to consider the question just in the isomorphic group $\mathrm{SU}(1,1)$ due to its nice structure, especially considering the fact that our main applications in this paper concern CMV matrices, with which one associates $\mathrm{SU}(1,1)$  cocycles. Let us summarize these results in $\mathrm{SU}(1,1)$; see Corollary \ref{openLinearly} below.

Suppose $\tilde{A}_{\Gamma}(x)\in \mathrm{SU}(1,1)$ and the cocycle $(\omega,M^{-1}\tilde{A}_{\Gamma}(x)M)$ obeys the hypotheses (H1)--(H3).
Let $\Gamma_{0}=(\theta_{0},\delta_{0})$ be such that there exists $Z(x)\in C^{\omega}(2\T^{d},\mathrm{SL}(2,\mathbb{R}))$ with
$$
Z^{-1}(x+\omega)M^{-1}\tilde{A}_{\Gamma_{0}}(x)MZ(x)=\begin{pmatrix}1 &c\\0 &1\end{pmatrix},
$$
which is equivalent to
$$
\tilde{Z}^{-1}(x+\omega)\tilde{A}_{\Gamma_{0}}(x)\tilde{Z}(x)=M\begin{pmatrix}1 &c\\0 &1\end{pmatrix}M^{-1},
$$
where $\tilde{Z} = M Z M^{-1}\in C^{\omega}(2\T^{d}, \mathrm{SU}(1,1))$.
For $\Gamma$ in a neighborhood of $\Gamma_{0}$, note that $\zeta=(\eta,\beta)=\Gamma-\Gamma_{0}$,
\begin{equation}\label{eq.linearPart}
\begin{aligned}\tilde{Z}^{-1}(x+\omega)\tilde{A}_{\Gamma}(x)\tilde{Z}(x)&=M\begin{pmatrix} 1&c\\0&1\end{pmatrix}M^{-1}\tilde{Z}^{-1}(x)\tilde{A}^{-1}_{\Gamma_{0}}(x)\tilde{A}_{\Gamma}(x)\tilde{Z}(x)\\
&=M\begin{pmatrix} 1&c\\0&1\end{pmatrix}M^{-1}e^{\tilde{P}(x,\zeta)+\tilde{O}_{2}(x,\zeta)},
\end{aligned}
\end{equation}
where $\tilde{P}(x,\zeta)$ denotes the linear part in $\zeta$ and $\tilde{O}_{2}(x,\zeta)$ collects the remainder of higher order in $\zeta$.

Then, as a consequence of Proposition~\ref{twoCases}, we have the following:

\begin{Corollary}\label{openLinearly}
Let $\omega \in \mathrm{DC}(\kappa,\tau), \Gamma_{0}=(\theta_{0},\delta_{0})$ be fixed. Assume that $(\omega,M^{-1}\tilde{A}_{\Gamma_{0}}(x)M)$ satisfies {\rm (H1)}--{\rm (H3)}, and assume that
$$
[\tilde{P}]=\begin{pmatrix} i [\tilde{P}_{11} ]& [\tilde{P}_{12} ] \\ [\overline{\tilde{P}_{12}}]&- i [\tilde{P}_{11}]\end{pmatrix} = \begin{pmatrix} i(a_1\eta+b_1\beta)&a_2\eta+b_2\beta \\\overline{a_2}\eta+\overline{b_2}\beta& -i(a_1\eta+b_1\beta) \end{pmatrix},
$$
where $a_1,b_1\in\mathbb{R},a_2,b_2\in\mathbb{C}$. Then we have the following:
\begin{enumerate}

\item  If $c\neq 0$, and $a_{1}\neq  \im \overline{a_{2}}$,  then the boundary $\theta=\theta(\delta)$ is real analytic in a neighborhood of $\delta_{0}$.

\item If $c=0$,   and if $a^{2}_{1}\neq |a_{2}|^{2}$,  then there exist two  tongue boundaries $\theta=\theta^{\pm}(\delta)$, which are real analytic in a neighborhood of $\delta_{0}$.

\item Furthermore, for $\delta_{0}=0$, if we assume $(a_1b_1-\Re a_2\overline{b_2})^{2}>(a_1^{2}-|a_2|^{2})(b_1^{2}-|b_2|^{2})$, then
\begin{equation}\label{eq.derivative}\frac{d\theta^{+}}{d\delta}(0)\neq\frac{d\theta^{-}}{d\delta}(0).
\end{equation}
In particular, if  $a_2=0$, then \eqref{eq.derivative} holds as long as $b_2\neq 0$.

\end{enumerate}
\end{Corollary}

\begin{pf}
Multiplying $M^{-1},M$ in the left and right side of \eqref{eq.linearPart}, we obtain
$$
\begin{aligned}
Z^{-1}(x+\omega)A_{\Gamma}Z(x)&=\begin{pmatrix}1 &c\\ 0 &1\end{pmatrix}Z^{-1}(x)A^{-1}_{\Gamma_{0}}(x)A_{\Gamma}(x)Z(x)\\
&=\begin{pmatrix}1 &c\\0 &1\end{pmatrix}e^{P(x,\zeta)+O_{2}(x,\zeta)}.
\end{aligned}
$$
It follows immediately that $[P(x,\zeta)]=[M^{-1}\tilde{P}(x,\zeta)M]=M^{-1}[\tilde{P}(x,\zeta)]M$,
and a simple calculation shows that
$$
[P]=\begin{pmatrix}\Re\overline{a_{2}}\eta+\Re\overline{b_{2}}\beta &(a_{1}+\im \overline{a_{2}})\eta+(b_{1}+\im\overline{b_{2}})\beta\\
(\im\overline{a_{2}} -a_{1} )\eta+  (\im\overline{b_{2}}-b_{1})\beta &-\Re\overline{a_{2}}\eta-\Re\overline{b_{2}}\beta
\end{pmatrix}.
$$
Therefore, we have
$$
[P_{21}]=(\im\overline{a_{2}} -a_{1} )\eta+  (\im\overline{b_{2}}-b_{1})\beta,
$$
$$
\det[P]=-(a_{1}^{2}-|a_{2}|^{2})\eta^{2}+O(\eta\beta,\beta^{2}).
$$
Then the statements (1) and (2) follow immediately from Proposition~\ref{twoCases}.

On the other hand, apply Theorem~\ref{KAM} to $(\omega, e^{P(x,\zeta)+O_{2}(x,\zeta)})$ to obtain a counterterm $\xi^{*}(\zeta)$, and averaging the cocycle $(\omega,e^{P(x,\zeta)+O_{2}(x,\zeta)}e^{-\xi^{*}(\zeta)})$ once, we get
$$
\xi^{*}=[P(x,\zeta)]+O_{2}(\zeta).
$$
Indeed, by the proof of Proposition~\ref{twoCases}, the two tongue boundaries $\theta=\theta^{\pm}(\delta)$ satisfy $\det\xi^{*}=0$, thus it is readily checked that their derivatives at $\beta=0$ are determined by the following equation,
$$
\det[P]=\det[\tilde{P}]=(a_1\eta+b_1\beta)^{2}-(a_2\eta+b_2\beta)(\overline{a_2}\eta+\overline{b_2}\beta)=0,
$$
which is equivalent to
$$
(a_1^{2}-|a_2|^{2})\eta^{2}+(2a_1b_1-2\Re(a_2\overline{b_2}))\eta\beta+(b_1^{2}-|b_2|^{2})\beta^{2}=0.
$$
Then \eqref{eq.derivative} follows immediately by considering its discriminant. This proves statement (3) and concludes the proof.
\end{pf}

\section{Cantor Spectrum for Quasi-Periodic CMV Matrices}\label{sec.4}

Let $\mathcal{E}_{\alpha}$ be the two-sided CMV matrix with quasi-periodic Verblunsky coefficients given by $\alpha_{n}= f (x+n\omega)$ where $f(x)=\lambda e^{ih(x)}$, $h\in C^{\omega}(\mathbb{T}^{d},\mathbb{R})$. We want to show for generic $f$ close to constant, the spectrum of $\mathcal{E}_{\alpha}$ has all gaps allowed by the gap labeling theorem open. To prove this, we can consider the corresponding \textit{Szeg\H{o} cocycles} $(\omega, S(\lambda e^{i \delta h(x)}, e^{i\theta}))$. Fix $\lambda$ and  $h$, and take $(\theta,\delta)$ as parameters. We first show that the tongue boundaries are analytic (see Proposition \ref{Thm.first}), and then proceed to show that the set of small $\delta$'s for which there are collapsed gaps is at most countable.

\subsection{Analytic Tongue Boundaries}

\begin{Proposition}\label{Thm.first}
Let $\lambda\in(0,1),\delta\in\mathbb{R},r>0$,  $\omega \in \mathrm{DC}(\kappa,\tau)$, $h\in C_{r}^{\omega}(\mathbb{T}^{d},\mathbb{R})$. Denote by $\mathcal{E}_{\alpha}$ the two-sided CMV matrix with quasi-periodic Verblunsky coefficients given by
$$
\alpha_{n}=\lambda e^{i \delta h(x+n\omega)}.
$$
Then there exists $\epsilon_{1} = \epsilon_{1}(\kappa,\tau,\|h\|_{r},\lambda) > 0$ such that if $|\delta_{0}|<\epsilon_{1}$ and the pair $(\theta_{0},\delta_{0})$ lies on a tongue boundary, we have the following:

(i) If $\theta_{0}$ is an endpoint of an open spectral gap of $\sigma(\mathcal{E}_{\alpha})$, then the tongue boundary $\theta=\theta(\delta)$ such that $\theta_{0}=\theta(\delta_{0})$ is real analytic in a neighborhood of $\delta_{0}$.

(ii) If $\theta_{0}$ lies at a collapsed spectral gap of $\sigma(\mathcal{E}_{\alpha})$, then the two tongue boundaries $\theta_{i}=\theta_{i}(\delta),i=1,2$ with $\theta_{i}(\delta_{0})=\theta_{0}$ are real analytic in a neighborhood of $\delta_{0}$.
\end{Proposition}

\begin{pf}
Since the Szeg\H{o} cocycle map belongs to $\mathrm{SU}(1,1)$, it is convenient to apply Corollary~\ref{openLinearly} instead of applying Proposition~\ref{twoCases} directly.

Note that $S(\lambda,e^{i\theta})$ is uniformly bounded in $\theta$, and hence we can take $|\delta|<\epsilon_1$, which is small enough such that
$$
\|S(\lambda e^{i\delta h},e^{i\theta})-S(\lambda,e^{i\theta})\|_{r}\leq \frac{\lambda}{\sqrt{1-\lambda^{2}}} \epsilon_1 \|h\|_{r}\leq  \epsilon_0(\kappa,\tau,r),
$$
where $\epsilon_0= \epsilon_0(\kappa,\tau,r)$ is defined in Theorem~\ref{Thm.reduci} (see also Remark~\ref{uniform}). Thus, (H3) holds.

Since $(\theta_{0},\delta_{0})$ lies on a tongue boundary, there exists $k\in\Z^d$ such that
$$
2 \rho(\omega,S(\lambda e^{i\delta_0 h},e^{i\theta_0}))=\langle k,\omega\rangle\mod\mathbb{Z}
$$
and $(\omega,S(\lambda e^{i\delta_0 h},e^{i\theta_0}))$ is not uniformly hyperbolic.

The non-decreasing property of the rotation number with respect to $\theta$ is due to a well known property of Szeg\H{o} cocycles \cite{Simon2}, and thus (H1) and (H2) hold as well. Therefore, it remains to verify conditions (1)--(2) of Corollary~\ref{openLinearly}.

Let $\alpha(x) = \lambda e^{i\delta h(x)}$ and
$$
S(\alpha,z)=\frac{z^{-\frac{1}{2}}}{\sqrt{1-\lambda^{2}}}\begin{pmatrix}z &-\lambda e^{-i\delta h(x)}\\
-\lambda e^{i\delta h(x)}z &1
\end{pmatrix}.
$$
Now we fix some $\delta=\delta_{0}, z_0=e^{i\theta_0}$ and assume that there exists
$$
\tilde{Z}(x) =  \begin{pmatrix} \tilde{z}_{11} &  \tilde{z}_{12} \\ \overline{\tilde{z}_{12}} & \overline{\tilde{z}_{11}} \end{pmatrix}      \in C^{\omega}(2\mathbb{T}^{d},\mathrm{SU}(1,1))
$$
such that
$$
\tilde{Z}^{-1}(x+\omega)S(\lambda e^{i\delta_{0}h},e^{i\theta_{0}})\tilde{Z}(x)=M\begin{pmatrix}1 & c\\ 0 &1\end{pmatrix}M^{-1}.
$$
Then for $(e^{i\theta},\delta)$ around $(e^{i\theta_0},\delta_{0})$, denote $\beta=\delta-\delta_{0}$,  $\eta=\theta-\theta_0$ and let $\zeta=(\eta,\beta)$ be the multi-parameter. Let $Q+O_{2}(x,\zeta)=\tilde{A}_{\Gamma_{0}}^{-1}(x)\tilde{A}_{\Gamma}(x)-Id$, where
$$
Q=\begin{pmatrix}i(\frac{\eta}{2}-\frac{\lambda^{2}}{1-\lambda^{2}}h \beta) &i\lambda e^{-i(\theta_{0}+\delta_{0}h)}\\ -i\lambda e^{i(\theta_{0}+\delta_{0}h)}h\beta &-i(\frac{\eta}{2}-\frac{\lambda^{2}}{1-\lambda^{2}}h\beta)\end{pmatrix},
$$
and let $\tilde{P}$ be as in \eqref{eq.linearPart}. Then, a direct computation shows that
$$
\tilde{P}=\tilde{Z}^{-1}Q\tilde{Z}=\begin{pmatrix}i(a_{1}(x)\eta+b_{1}(x)\beta) &a_{2}(x)\eta+b_{2}(x)\beta\\
\overline{a_{2}(x)}\eta+\overline{b_{2}(x)}\beta &-i(a_{1}(x)\eta+b_{1}(x)\beta) \end{pmatrix},
$$
where
\begin{equation}\label{eq.coefficients1}
\begin{aligned}
&a_{1}(x)=\frac{|\tilde{z}_{11}|^{2}+|\tilde{z}_{12}|^{2}}{2}\eta\\
&b_{1}(x)=\frac{\lambda h}{1-\lambda^{2}}(2\Re(\tilde{z}_{11}\tilde{z}_{12}e^{i(\theta_{0}+\delta_{0}h)})-\lambda(|\tilde{z}_{11}|^{2}+|\tilde{z}_{12}|^{2}))\\
&a_{2}(x)=\overline{\tilde{z}_{11}}\tilde{z}_{12}\\
&b_{2}(x)=\frac{\lambda h}{1-\lambda^{2}}(\tilde{z}_{12}^{2}e^{i(\theta_{0}+\delta_{0}h)}+\overline{\tilde{z}_{11}}^{2}e^{-i(\theta_{0}+\delta_{0}h)}-2\lambda\overline{\tilde{z}_{11}}\tilde{z}_{12})
\end{aligned}
\end{equation}

We obtain
$$
a_{1}=[a_{1}(x)]=\frac{[|\tilde{z}_{11}|^{2}+|\tilde{z}_{12}|^{2}]}{2},\quad  a_{2}=[a_{2}(x)]=[\overline{\tilde{z}_{11}}\tilde{z}_{12}].
$$

We are now ready to verify the assumptions of Corollary~\ref{openLinearly}: Since $|\tilde{z}_{11}|^{2}-|\tilde{z}_{12}|^{2}=1$, we always have $a_{1}=\frac{[|\tilde{z}_{11}|^{2}+|\tilde{z}_{12}|^{2}]}{2}>|a_{2}|=|[\overline{\tilde{z}_{11}}\tilde{z}_{12}]|$, which immediately verifies both (1) and (2) in Corollary~\ref{openLinearly}. This completes the proof of Proposition \ref{Thm.first}.
\end{pf}

\subsection{Cantor Spectrum}\label{Cantor}

We first note that in the case $\delta=0$, the sequence of Verblunsky coefficients is constant. It is well known that in this situation, the spectrum is given explicitly by $\{e^{i\theta}:2\arcsin |\lambda|\leq\theta\leq 2\pi-2\arcsin |\lambda|\}$. Then $e^{i\theta_{0}}$ is clearly an edge of the unique open gap with label $k=0$, since the rotation number equals $0$ in this gap. The gap of course remains open under small perturbations. In the following argument, we exclude this case from our considerations and assume that $k\neq 0$.

In order to prove generic Cantor spectrum with all gaps open, we need to show that the two tongue boundaries of a given resonance have a transversality at the origin, which is the content of the following result:

\begin{Proposition}\label{openGap}
Denote the tongue boundaries of the gap with label $k\neq 0$ by $\theta_{k}^{\pm}(\delta)$. Then we have $\frac{d\theta^{+}_{k}}{d\delta}(0)\neq \frac{d\theta^{-}_{k}}{d\delta}(0)$ if and only if
$\hat{h}_{k}\neq 0$.
\end{Proposition}

\begin{pf}
Let $\delta_{0}=0, z=e^{i\theta_{0}}$. Then the matrix becomes
$$
S_{0}=\frac{1}{\sqrt{1-\lambda^{2}}}\begin{pmatrix}z^{\frac{1}{2}} &-\lambda z^{-\frac{1}{2}}\\-\lambda z^{\frac{1}{2}} &z^{-\frac{1}{2}}\end{pmatrix}.
$$
As we are considering tongue boundaries of the gap with label $k$, this means that the eigenvalues of $S_0$ must be given by $\lambda_{\pm}(\theta)=e^{\pm i\frac{\langle k,\omega\rangle}{2}}$. Indeed, by direct computation one can show that the two eigenvalues are given by
$$
\lambda_{\pm}(\theta)=\frac{\cos\frac{\theta_{0}}{2}}{\sqrt{1-\lambda^{2}}}\pm i\sqrt{1-(\frac{\cos\frac{\theta_{0}}{2}}{\sqrt{1-\lambda^{2}}})^{2}},
$$
and one can diagonalize $S_0$ with the help of the following result:

\begin{Lemma}[\cite{KXZ}]\label{diagonalizing}
Let $A = \begin{pmatrix}it &z\\\bar{z} &-it\end{pmatrix}\in \mathrm{su}(1,1)$ with $t\in\mathbb{R},z\in\mathbb{C}$. Assume that $\det A>0$ and let $\rho=\sqrt{\det A}$. Then we have
$$
U^{-1}AU=\begin{pmatrix} i\rho & 0 \\ 0 & -i\rho\end{pmatrix},
$$
where $U$ takes the form
\begin{eqnarray}\label{formu}
U = \frac{1}{\sqrt{\cos2 \tilde{\theta}}}\begin{pmatrix}\cos\tilde{\theta} &e^{2i\phi}\sin \tilde{\theta}\\e^{-2i\phi}\sin \tilde{\theta} &\cos\tilde{\theta}\end{pmatrix}.
\end{eqnarray}
Here, $2\phi=\arg z-\frac{\pi}{2}$ and $\tilde{\theta}\in(-\frac{\pi}{2},\frac{\pi}{2})$ satisfies
$$
2\tilde{\theta}=-\arctan\frac{|z|}{\sqrt{t^{2}-|z|^{2}}}.
$$
In addition we have
$$
\|U\|^{2}= \frac{|t|+|z|}{\rho}.
$$
\end{Lemma}

Let $U$ be given by Lemma \ref{diagonalizing} so that
$$
U^{-1}S_{0}U=\begin{pmatrix}e^{i\frac{\langle k,\omega\rangle}{2}} &0\\0 &e^{-i\frac{\langle k,\omega\rangle}{2}}\end{pmatrix}.
$$
Let $\tilde{Z}(x)=U\begin{pmatrix}e^{i\frac{\langle k,x\rangle}{2}} &0\\0 &e^{-i\frac{\langle k,x\rangle}{2}}\end{pmatrix}.$ Then,
$$
\tilde{Z}^{-1}(x+\omega)S_{0}\tilde{Z}(x)=Id,
$$
and it conjugates the perturbed cocycle $(\omega,S(\lambda e^{i\beta h},ze^{i\eta}))$ into
$$
\tilde{Z}^{-1}(x+\omega)   S(\lambda e^{i\beta h},ze^{i\eta}) \tilde{Z}(x)=e^{\tilde{P}(x,\zeta)+\tilde{O}_{2}(x,\zeta)}.
$$
Now we only need to check the elements of $\tilde{P}$.

Note that since $U$ takes the form \eqref{formu}, we have \begin{equation}\label{eq42201}
\tilde{Z}(x)=\frac{1}{\sqrt{\cos2\tilde{\theta}}}\begin{pmatrix}\cos\tilde{\theta}e^{i\frac{\langle k,x\rangle}{2}} &e^{2i\phi}\sin\tilde{\theta}e^{-i\frac{\langle k,x\rangle}{2}}\\
e^{-2i\phi}\sin\tilde{\theta}e^{i\frac{\langle k,x\rangle}{2}} &\cos\tilde{\theta}e^{-i\frac{\langle k,x\rangle}{2}}
\end{pmatrix}.
\end{equation}
According to \eqref{eq.coefficients1} with $\delta_{0}=0,z=e^{i\theta_{0}}$ and the elements of $\tilde{Z}$ above, we have
$$
\begin{aligned}a_2(x)&= \frac{i e^{2i\phi}\cos\tilde{\theta}\sin\tilde{\theta}e^{-i\langle k,x\rangle}}{\cos2\tilde{\theta}   }  \\
         b_2(x)&=\frac{ -i \lambda^{2}e^{2i\phi}\sin2\tilde{\theta} + i\lambda (e^{i(\theta_{0}+4\phi)}\sin^{2}\tilde{\theta}+e^{-i\theta_{0}}\cos^{2}\tilde{\theta})}{(1-\lambda^{2})\cos2\tilde{\theta}   }e^{-i\langle k,x\rangle}h.
\end{aligned}
$$

Taking the average, we obtain $a_2=[a_2(x)]=0$ and
$$
b_2= [b_2(x)]=i\frac{\lambda}{1-\lambda^{2}}((e^{i(\theta_{0}+4\phi)}\sin^{2}\tilde{\theta}+e^{-i\theta_{0}}\cos^{2}\tilde{\theta})-\lambda e^{2i\phi}\sin2\tilde{\theta})\hat{h}_{-k}.
$$
Thus by Corollary \ref{openLinearly} (2),(3), there exist two analytic tongue boundaries $\theta_{k}^{\pm}(\delta)$ and $\frac{d\theta^{+}_{k}}{d\delta}(0)\neq \frac{d\theta^{-}_{k}}{d\delta}(0)$ if and only if $b_2 \neq 0$. Just note
\begin{eqnarray*}|b_{2}| &=&\frac{\lambda|\hat{h}_{k}|}{(1-\lambda^{2}) |\cos2\tilde{\theta}|}|(e^{i(\theta_{0}+4\phi)}\sin^{2}\tilde{\theta}+e^{-i\theta_{0}}\cos^{2}\tilde{\theta})-\lambda e^{2i\phi}\sin2\tilde{\theta}| \\
&=&\frac{\lambda|\hat{h}_{k}\cos\tilde{\theta}|}{(1-\lambda^{2}) |\cos2\tilde{\theta}|}|(e^{i(\theta_{0}+2\phi)}\tan\tilde{\theta}-\lambda)^{2}+1-\lambda^{2}|,
 \end{eqnarray*}
and according to Lemma \ref{diagonalizing}, $\tilde{\theta}\in(-\frac{\pi}{2},\frac{\pi}{2})$ never equals $\frac{\pi}{4}$. Hence $b_2\neq 0$ if and only if $\hat{h}_{k}\neq 0$, finishing the proof.
\end{pf}

\medskip
\noindent\textit{Proof of  Theorem \ref{MainThm.refined}.} Let
$$
\mathcal{G}_{k} = \left\{ h \in C^{\omega}(\mathbb{T}^{d},\mathbb{R}) : \hat{h}_{k}\neq 0 \right\}
$$
and
$$
\mathcal{G} = \bigcap_{k\in\mathbb{Z}^{d}}\mathcal{G}_{k}.
$$
Clearly, $\mathcal{G}$ is a dense $G_\delta$ set.

For any $h\in \mathcal{G}$ and any label $k \not= 0$, the set of couplings $|\delta|<\epsilon$ for which the tongue boundaries corresponding to label $k$ coincide is finite, since the tongue boundary functions are real analytic functions of $\delta$ (see Proposition~\ref{Thm.first}), and the tongue boundaries  have a transversality at the origin (see Proposition~\ref{openGap}).  This completes the proof of Theorem \ref{MainThm.refined}.
\qed

\section{Cantor Spectrum for Almost Periodic CMV Matrices}\label{AlmostModel}

In the previous section, we proved the analyticity of the tongue boundaries and generic Cantor spectrum for a class of two-sided quasi-periodic CMV matrices. In this section, we work out results of this type for a class of almost periodic CMV matrices.

\subsection{The Model and its Gap Labelling}

Consider the  torus $\mathbb{T}^{d}$ equipped with Lebesgue measure and $\mathbb{Z}_{2}$ equipped with the probability measure that assigns the weight $\frac12$ to each of $0$ and $1$. Let $\Omega=\mathbb{T}^{d}\times \mathbb{Z}_{2}$ be the product space and $\nu$ the product measure. The transformation $T_{\omega}: \Omega\to\Omega$ is given by $(x,j)\to (x+\omega,j+1)$, where $\omega$ is assumed to be Diophantine. It is readily verified that $T_{\omega}$ is uniquely ergodic with unique invariant  measure $\nu$ (compare, e.g., Theorem 9.1 of \cite{mane2012ergodic}).

Given $h\in C^{\omega}(\mathbb{T}^{d},\mathbb{R})$, $\lambda_1, \lambda_2 \in [0,1)$ {\color{blue}with $\lambda_{1}^{2}+\lambda_{2}^{2}\neq 0$,} and $\delta \in \mathbb{R}$, we consider the Verblunsky coefficients given by
\begin{equation}\label{eq.sampling-1}
f(x,j)=\left\{\begin{aligned} \lambda_{1}e^{ i\delta h(x)}\qquad \text{for } j=0 \mod 2,\\\lambda_{2}e^{i \delta h(x)}\qquad \text{for } j=1 \mod 2.\end{aligned}\right.
\end{equation}
and \begin{equation}\label{eq.way-1}
\alpha_{n}(x) = f(x + n\omega, n), \; x \in \T^d , \; n \in \mathbb{Z},
\end{equation}
 and write $\alpha=\{\alpha_{n}\}_{n\in\mathbb{Z}}$ for short. The associated cocycle is denoted by $(T_\omega,S(\alpha,z))$, where $S(\alpha,z)$ is the usual Szeg\H{o} cocycle map and $z$ is the spectral parameter.

We want to show that for generic $f$ close to constant,  the corresponding CMV matrix has all gaps open. To make this more precise, let us recall what the Gap Labelling Theorem says about the gap labels that occur. While this could be done based on Johnson-Moser's approach \cite{JM82}, we will give a self-contained proof based on methods from dynamical systems. The useful observation is that while $(T_\omega, S(\alpha,z))$ is not a quasi-periodic cocycle, its iterate
$$
(2\omega, T_{\delta}(\lambda_2,z)(\cdot))=: (2\omega, S(\alpha_1(\cdot),z)S(\alpha_0(\cdot),z))
$$
indeed does define an analytic quasi-periodic cocycle. Here, we make the dependence of $T_{\delta}(\lambda_2,z)$ on $\lambda_2$ explicit, since this dependence will be quite important in the proof.  As a result of this, we can describe the labelling associated with $(T_\omega,S(\alpha,z))$ with the help of the labelling associated with suitable quasi-periodic $\SL(2,\R)$ cocycles.

\begin{Proposition}\label{Lem.gapLabel}
Let $\omega\in \T^d$ be rationally independent and let $\mathcal{I}\subset\partial\mathbb{D}$ be an open arc. Then for any  $e^{i\theta}\in\mathcal{I}$, if $(T_\omega,S(\alpha,e^{i\theta}))$ is uniformly hyperbolic, there exists $k\in \mathbb{Z}^{d}$ such that either
\begin{equation}\label{labeling1}
2\rho(T_\omega,S(\alpha,e^{i\theta}))=\langle k,\omega\rangle \mod \mathbb{Z}
\end{equation}
or
\begin{equation}\label{labeling2}
2\rho(T_\omega,S(\alpha,e^{i\theta}))=\langle k,\omega\rangle +\frac{1}{2} \mod \mathbb{Z}
\end{equation}
\end{Proposition}
\begin{pf}
First we need the following simple observation:
   \begin{Lemma}\label{eqsystem}\cite{V2}
For any $z\in\partial \mathbb{D}$, the system $(T_\omega,S(\alpha,e^{i\theta}))$ is uniformly hyperbolic if and only if the system $(2\omega,T_{\delta}(\lambda_2,e^{i\theta}))$ is uniformly hyperbolic.
\end{Lemma}
Therefore, under the assumption of the proposition, we have that the cocycle $(2\omega, T_{\delta}(\lambda_2,e^{i\theta}))$ is uniformly hyperbolic for any $e^{i\theta}\in \mathcal{I}$. Thus there exists an analytic $B:2\T^{d}\to \mathrm{SU}(1,1)$ such that
$$
B^{-1}(x+2\omega)  T_{\delta}(\lambda_2,e^{i\theta})(x) B(x)= \begin{pmatrix}
  \tilde{\lambda}(x) &0\\
  0 &   \tilde{\lambda}^{-1}(x)
  \end{pmatrix} .
$$
Assume that the degree of $B$ is $k\in \mathbb{Z}^{d}$. Then,
$$
2\rho(2\omega,T_{\delta}(\lambda_{2},e^{i\theta}))=\langle k,2\omega\rangle \mod\mathbb{Z}.
$$

Note that the rotation numbers of the two cocycles $(T_\omega, S(\alpha, e^{i\theta}))$ and $(2\omega, T_{\delta}(\lambda_{2}, e^{i\theta}))$ satisfy
\begin{equation}\label{double}
\rho(2\omega,T_{\delta}(\lambda_{2}, e^{i\theta})) = 2 \rho(T_\omega, S(\alpha, e^{i\theta})) \mod \mathbb{Z}.
\end{equation}
It follows immediately that either
$$
2 \rho(T_\omega,S(\alpha,e^{i\theta})) = \langle k, \omega \rangle \mod \mathbb{Z}
$$
or
$$
2 \rho(T_\omega, S(\alpha,e^{i\theta})) = \langle k, \omega \rangle + \frac{1}{2} \mod \mathbb{Z}.
$$
This finishes the proof.
\end{pf}

\subsection{Analytic Tongue Boundaries}\label{period2}

Let us now prove that the tongue boundary of $(2\omega,T_{\delta}(\lambda_2,e^{i\theta}))$ is analytic.

\begin{Proposition}\label{Thm.period2}
Let $\omega \in \mathrm{DC}(\kappa,\tau)$ and $h\in C^{\omega}(\mathbb{T}^{d},\mathbb{R})$.  Suppose that $\alpha$ is given by \eqref{eq.way-1}.
Then, for each fixed $\lambda_{1}\in[0,1)$, there exist $r_{1}(\lambda_{1})\leq \frac{1}{2}$ and $\tilde{\epsilon}= \tilde{\epsilon}(\kappa,\tau,h,\lambda_{1}) > 0$ such that for any $|\lambda_{2} | \leq r_{1}(\lambda_{1})$, $|\delta|< \tilde{\epsilon},$  if $(\theta_{0},\delta_{0})\in\mathbb{R}^{2}$ lies on a tongue boundary of  $(2\omega,T_{\delta}(\lambda_2,e^{i \theta}))$, then we have the following:

(i) If the tongue is non-degenerate,   then the boundary $\theta(\delta)$ such that $\theta(\delta_{0})=\theta_{0}$ is real analytic in a neighborhood of $\delta_{0}$.

(ii) If the tongue is degenerate, then the two tongue boundaries $\theta_{i}=\theta_{i}(\delta)$ for $i=1,2$ with $\theta_{i}(\delta_{0})=\theta_{0}$ are real analytic in a neighborhood of $\delta_{0}$.
\end{Proposition}

\begin{pf}
The strategy of the proof is to apply Corollary~\ref{openLinearly} to the  iterated quasi-periodic cocycle $(2\omega,T_{\delta}(\lambda_2,e^{i\theta}))$. Note that if $\omega \in \mathrm{DC}(\kappa,\tau)$, then $2\omega \in \mathrm{DC}( \frac{\kappa}{2^{\tau}},\tau)$.

We first consider the case where $\lambda_{1}\in (0,1)$ is fixed and $\lambda_2 \in [0,\frac{1}{2}]$. A direct computation gives \begin{equation}\label{twosteps}
\begin{aligned}
&T_{\delta}(\lambda_{2},e^{i\theta})=\\&\frac{1}{\sqrt{(1-\lambda_{1}^{2})(1-\lambda_{2}^{2})}}\begin{pmatrix}e^{i\theta}+\lambda_{1}\lambda_{2}e^{i\delta(h_{+}-h)} &-\lambda_{1}e^{-i(\theta+\delta h)}-\lambda_{2}e^{-i\delta h_{+}}\\
-\lambda_{1}e^{i(\theta+\delta h)}-\lambda_{2}e^{i\delta h_{+}} &e^{-i\theta}+\lambda_{1}\lambda_{2}e^{-i\delta(h_{+}-h)}
\end{pmatrix},
\end{aligned}
\end{equation}
where we simply denote $h_{+}=h(x+\omega)$. Note that $T_{0}(\lambda_{2},e^{i\theta})$ is uniformly bounded in $\theta$. Thus we can take $|\delta|<\epsilon_2$ which is small enough so that
$$
\|T_{\delta}(\lambda_{2},e^{i\theta})-T_{0}(\lambda_{2},e^{i\theta})\|_{r}\leq \frac{1}{\sqrt{(1-\lambda_{1}^{2})(1-\lambda_{2}^{2})}} |\delta| \|h\|_{r} \leq   \epsilon_0( \frac{\kappa}{2^{\tau}},\tau,r),
$$
where $\epsilon_0= \epsilon_0(\kappa,\tau,r)$ is defined in Theorem \ref{Thm.reduci} (see also Remark~\ref{uniform}). Thus, (H3) holds.
Since $(\theta_{0},\delta_{0})$ lies on a tongue boundary, which means that there exists $k\in\Z^d$ such that
$$
2\rho(2\omega,T_{\delta_0}(\lambda_{2},e^{i\theta_0}))=\langle k,2\omega\rangle \mod\mathbb{Z}
$$
and $(2\omega,T_{\delta_0}(\lambda_{2},e^{i\theta_0}))$ is not uniformly hyperbolic. On the other hand, by \eqref{double} and  $\rho(T_\omega, S(\alpha,e^{i\theta}))$ being non-decreasing with respect to $\theta$ \cite{Simon2}, we know that  $\rho(2\omega,T_{\delta}(\lambda_{2},e^{i\theta}))$ is not locally monotonic if and only if $$\rho(T_{\omega}, S(\alpha,e^{i\theta})) = \frac{1}{4}.$$ This directly implies that $\rho(2\omega,T_{\delta_0}(\lambda_{2},e^{i\theta_0}))$ is a non-decreasing (or non-increasing) function of $\theta$ in a neighborhood of $\theta_{0}$. Thus (H1) and (H2) hold as well, and it remains to verify conditions (1) and (2) of Corollary~\ref{openLinearly}.

Let $\Gamma_{0}=(\theta_{0},\delta_{0})$ and denote $\tilde{A}_{\Gamma}(x)=T_{\delta}(\lambda_2,e^{i\theta})$. Then there exists $\tilde{Z}\in C^{\omega}(2\mathbb{T}^{d},\mathrm{SU}(1,1))$ such that
\begin{equation}\label{eq.firstCon}
\tilde{Z}^{-1}(x+\tilde{\omega})\tilde{A}_{\Gamma_{0}}(x)\tilde{Z}(x)=M\begin{pmatrix}1 &c\\0 &1\end{pmatrix}M^{-1}.
\end{equation}
For $\Gamma$ in a neighborhood of $\Gamma_{0}$, we still denote $\zeta=(\eta,\beta)=\Gamma-\Gamma_{0}$. A direct computation shows that
\begin{eqnarray*}
\tilde{A}_{\Gamma_{0}}^{-1}(x)\tilde{A}_{\Gamma}(x)-Id  =  \begin{pmatrix}i Q_{11} & Q_{12}\\\overline{Q_{12}} &-i Q_{11}\end{pmatrix}+O_2(x,\zeta),
\end{eqnarray*}
where
\begin{equation}\label{eq.suCal4}\begin{aligned}
&Q_{11}=\frac{(1-\lambda_{1}^{2})\eta+f_{11}\beta }{(1-\lambda_{1}^{2})(1-\lambda_{2}^{2})},\\
& Q_{12}=\frac{i \lambda_{2}e^{-i(\theta_{0}+\delta_{0}h_{+})}(\lambda_{1}^{2}-1)\eta+ i f_{12}\beta }{(1-\lambda_{1}^{2})(1-\lambda_{2}^{2})}.
\end{aligned}
\end{equation}
Here $f_{11}$ and $f_{12}$ are defined by
\begin{equation}\label{eq.coefficients3}
f_{11}=\lambda_{1}^{2}\lambda_{2}^{2}(h_{+}-h)-\lambda_{1}^{2}h-\lambda_{2}^{2}h_{+}-2\lambda_{1}\lambda_{2}h\cos(\theta_{0}-\delta_{0}(h_{+}-h))
\end{equation}
and
\begin{equation}\label{eq.coefficients4}
\begin{aligned}
f_{12}&=(\lambda_{1}e^{-i(2\theta_{0}+\delta_{0}h)}+\lambda_{1}^{2}\lambda_{2}e^{-i(\theta_{0}+\delta_{0}h_{+})})h\\
&+\lambda_{2}e^{-i\delta_{0}h_{+}}(e^{-i\theta_{0}}+\lambda_{1}\lambda_{2}e^{-i\delta_{0}(h_{+}-h)})h_{+}\\
&-\lambda_{1}\lambda_{2}e^{-i\delta_{0}(h_{+}-h)}(\lambda_{1}e^{-i(\theta_{0}+\delta_{0}h)}+\lambda_{2}e^{-i\delta_{0}h_{+}})(h_{+}-h).
\end{aligned}
\end{equation}

Consequently, we have
\begin{equation*}
\begin{aligned}
\tilde{Z}^{-1}(x+\omega)\tilde{A}_{\Gamma}(x)\tilde{Z}(x)
=M\begin{pmatrix} 1&c\\0&1\end{pmatrix}M^{-1}e^{\tilde{P}(x,\zeta)+\tilde{O}_{2}(x,\zeta)}
\end{aligned}.
\end{equation*}
Indeed, a direct computation shows that
\begin{eqnarray*}
\tilde{P}=\tilde{Z}^{-1}Q\tilde{Z} =  \begin{pmatrix} i a_1 (\lambda_2,x)\eta+ i b_1(\lambda_2,x)\beta &a_2 (\lambda_2,x)\eta+b_2(\lambda_2,x)\beta \\\overline{a_2(\lambda_2,x)}\eta+\overline{b_2(\lambda_2,x)}\beta& - i(a_1(\lambda_2,x)\eta-i b_1( \lambda_2,x)\beta \end{pmatrix} ,
\end{eqnarray*}
where
 \begin{equation}\label{Peq1}\begin{aligned}
a_1(\lambda_2,x)&=\frac{|\tilde{z}_{11}|^{2}+|\tilde{z}_{12}|^{2}}{1-\lambda_{2}^{2}}-\frac{\lambda_{2}}{1-\lambda_{2}^{2}}2\Re(\tilde{z}_{11}\tilde{z}_{12}e^{i(\theta_{0}+\delta_{0}h_{+})}),\\
b_1(\lambda_2,x) &= \frac{1}{(1-\lambda_{1}^{2})(1-\lambda_{2}^{2})}(f_{11}(|\tilde{z}_{11}|^{2}+|\tilde{z}_{12}|^{2})+2\Re(\overline{\tilde{z}_{11}\tilde{z}_{12}}f_{12})),\\
a_2(\lambda_2,x)&=\frac{1}{1-\lambda_{2}^{2}}\left(
-\lambda_{2}(\overline{\tilde{z}_{11}}^{2}e^{-i(\theta_{0}+\delta_{0}h_{+})}+\tilde{z}_{12}^{2}e^{i(\theta_{0}+\delta_{0}h_{+})})+2\overline{\tilde{z}_{11}}\tilde{z}_{12} \right),\\
b_2(\lambda_2,x)&=\frac{i}{(1-\lambda_{1}^{2})(1-\lambda_{2}^{2})}\left(2\overline{\tilde{z}_{11}}\tilde{z}_{12}f_{11}+\tilde{z}_{12}^{2}\overline{f_{12}}+\overline{\tilde{z}_{11}}^{2}f_{12}\right).
\end{aligned}
\end{equation}
Taking averages, we obtain
$$
a_{1}(\lambda_2)=[a_1(\lambda_2,x)]= \frac{1}{1-\lambda_{2}^{2}}[|\tilde{z}_{11}|^{2}+|\tilde{z}_{12}|^{2}-2\lambda_{2}\Re(\tilde{z}_{11}\tilde{z}_{12}e^{i(\theta_{0}+\delta_{0}h_{+})})],
$$
$$
a_{2}(\lambda_2)=[a_2(\lambda_2,x)]=\frac{1}{1-\lambda_{2}^{2}}[-\lambda_{2}(\overline{\tilde{z}_{11}}^{2}e^{-i(\theta_{0}+\delta_{0}h_{+})}+\tilde{z}_{12}^{2}e^{i(\theta_{0}+\delta_{0}h_{+})})+2\overline{\tilde{z}_{11}}\tilde{z}_{12}].
$$

First note that by Schwarz's inequality, we always have $$[|\tilde{z}_{11}|^{2}+|\tilde{z}_{12}|^{2}]>2|[\overline{\tilde{z}_{11}}\tilde{z}_{12}]|,$$ since $|\tilde{z}_{11}|^{2}-|\tilde{z}_{12}|^{2}=1$, which implies $a_{1}(0)>|a_{2}(0)|$. In order to prove
$$
a_{1}(\lambda_2)>|a_{2}(\lambda_2)|
$$
for $\lambda_{2}$ in a neighborhood of $0$, we just need to show that the $\tilde{z}_{i,j}$ are continuous in $\lambda_{2}$.

Note that for any fixed $\lambda_1, \delta$, and fixed $\lambda_2$,  there exists $e^{i\theta_{0}}\in \partial \mathbb{D}$, such that
$$
2\rho(2\omega,T_{\delta}(\lambda_2,e^{i\theta_{0}}))=\langle k,2\omega\rangle \mod\mathbb{Z}.
$$
Eliasson's result (Theorem \ref{Thm.reduci}) ensures that $(2\omega,T_{\delta}(\alpha,e^{i\theta}) $ is reducible by the conjugation $Z$. If we only perturb $\lambda_2$, clearly $Z$ is not well-defined (since if the rational number changes, once it becomes Liouvillean, the cocycle is not reducible any more). However, we will show that if we fix $\lambda_1,\delta$, and perturb $(\lambda_2,e^{i\theta})$ such that
$$
2\rho(2\omega,T_{\delta}(\lambda_2,e^{i\theta}))=\langle k,2\omega\rangle \mod\mathbb{Z},
$$
$Z$ is actually continuous in $(\lambda_2,\theta)$.

Denote $Z=Z(\mu,x)$, where $\mu=(\lambda_{2},\theta)$, to make this dependence explicit and denote $T_{\delta}(\mu,x)=T_{\delta}(\lambda_{2},e^{i\theta})$. For any $\rho>0,$ define
\begin{equation}\label{eq.paraSet}
\mathcal{R}_{\rho} = \left\{\mu\in[0, \tfrac{1}{2}]\times\partial\mathbb{D}:\rho(2\omega, T_{\delta}(\mu,x))=\rho \right\}.
\end{equation}
Then we have the following:

\begin{Proposition}\label{conDep}
Let $\lambda_1\in (0,1)$, $r>0$, and $\lambda_2\in[0,\frac{1}{2}]$.   Assume that $\omega \in \mathrm{DC}(\kappa,\tau)$,   $\tilde{\rho}$ is Diophantine w.r.t $2\omega$, or
$$
2\tilde{\rho}=\langle k,2\omega\rangle\mod\mathbb{Z}.
$$
Then there exists $\epsilon_{3} = \epsilon_{3}(\kappa, \tau, r,\lambda_{1}) > 0$, which is independent of $\lambda_2$, such that if $|\delta|\leq \epsilon_{3}$, there exist $Z(\mu,x)\in C^{0}(\mathcal{R}_{\tilde{\rho}}\times 2\T^{d},\mathrm{SU}(1,1))$, $C(\mu) \in C^0(  \mathcal{R}_{\tilde{\rho}}, \mathrm{SU}(1,1))$  such that
$$
Z^{-1}(\mu,x+2\omega)T_{\delta}(\mu,x)Z(\mu,x)=C(\mu).
$$
\end{Proposition}

We just take
$$
|\delta|\leq  \min\{  \epsilon_{2}(\kappa,\tau,h,\lambda_{1}), \epsilon_{3}(\kappa,\tau,h,\lambda_{1})\}.
$$
Then, by Proposition \ref{conDep}, $Z(\mu,x)$ is continuous in the set $\mathcal{R}_{\tilde{\rho}}$. Thus if $| \lambda_2|\leq r_{1}(\lambda_1)$, which is small enough, then $a_{1}(\lambda_2)>|a_{2}(\lambda_2)|$. Just note that $r_{1}(\lambda_1)$ can be chosen to be independent of $\delta$, the key observation here is that $\delta$ varies in a compact set, and this compact set is independent of $\lambda_2$; the result then follows from simple compactness argument.

The condition $a_{1}(\lambda_2)>|a_{2}(\lambda_2)|$ immediately implies condition (1) and (2) of Corollary~\ref{openLinearly}, finishing the proof.
\end{pf}

\begin{Remark}
Proposition~\ref{conDep} will be proved by KAM iteration, the key iteration step is Lemma~\ref{eq.oneKam} below. We note that Eliasson's scheme \cite{Eli92} was based on a parameterized KAM scheme, where the parameter is usually the energy $E$, thus the rotation number already has some monotonicity with respect to the parameter, which is quite different from our situation. Our scheme here (inspired by \cite{CCYZ,LYZZ}) is different from Eliasson's scheme \cite{Eli92}, where we can just fix the rotation number, and one can easily see how the conjugation $Z(\mu,x)$ varies when $\mu$ varies in $\mathcal{R}_{\rho}$.
\end{Remark}

\subsection{Proof of Proposition \ref{conDep}}

Before giving the proof, we first introduce some useful notations. Let $\mathcal{O}$ be a compact set of external parameters $\mu$, and for any  $r>0$, we denote by $\mathcal{B}_{\mathcal{O},r}(*)$ the set of $(*)$-valued functions\footnote{Here $(*)$ may denote a group such as $\mathrm{SU}(1,1),\SL(2,\mathbb{R})$ or an algebra such as $\mathrm{su}(1,1),\sl(2,\mathbb{R})$.}  $f(\mu,\cdot)\in C^{0}(\mathcal{O}\times\T^{d},*)$ and for each $\mu\in\mathcal{O}$, $f(\mu,\cdot)\in C^{\omega}_{r}(2\mathbb{T}^{d},*)$.

Our iteration step can now be stated as follows:

\begin{Lemma}\label{eq.oneKam}
Let $\tilde{\omega}\in DC(\tilde{\kappa},\tau)$, $r>0$, $A(\mu)\in C^{0}(\mathcal{O},\mathrm{SU}(1,1))$, and $f(\mu,\cdot)\in \mathcal{B}_{\mathcal{O},r}(\mathrm{su(1,1)})$. Assume that
$$
\rho(\tilde{\omega}, A(\mu)e^{ f(\mu,\cdot)} ) =\rho, \quad \text{for any} \quad  \mu \in \mathcal{O}.
$$
Then for any $r'<r$, there exist a numerical constant $C_{0}$ and a constant $D_{0}=D_{0}( \tilde{\kappa},\tau,d)$ such that if \begin{equation}\label{initialcon}
\|f(\mu,\cdot)\|_{r}\leq\epsilon\leq \frac{D_{0}}{\sup_{\mu}\|A(\mu)\|^{C_{0}}}(\min\{1,\frac{1}{r}\}(r-r'))^{C_{0}\tau},
\end{equation}
there exists $Z(\mu,\cdot)\in \mathcal{B}_{\mathcal{O},r'}(\mathrm{SU}(1,1))$ such that
$$
Z^{-1}(\mu,x+\tilde{\omega})A(\mu)e^{f(\mu,x)}Z(\mu,x)=A_{+}(\mu)e^{f_{+}(\mu,x)}
$$
with $\|f_{+}\|_{r'}\leq \epsilon^{2}$. More precisely, letting $N = \frac{2}{r-r'}|\ln\epsilon|$, we can distinguish two cases:
\begin{enumerate}

\item (Non-resonant case) If for any $n \in \mathbb{Z}^{d}$ with $0<|n|\leq N$, we have $$\|2\rho-\langle n,\tilde{\omega}\rangle\|_{\mathbb{R}/\mathbb{Z}}\geq\epsilon^{\frac{1}{15}},$$
then we have $\|Z(\mu,\cdot)-Id\|_{r'}\leq \epsilon^{\frac{1}{2}} ,$  $\lVert A_+ (\mu)-A(\mu) \rVert\leqslant 2\lVert A(\mu)\rVert\epsilon,$  and
$$
\rho(\tilde{\omega}, A_+(\mu)e^{ f_+(\mu,\cdot)} ) =\rho, \quad \text{for any} \quad  \mu \in \mathcal{O}.
$$

\item (Resonant case) If for some $n_{*}\in\mathbb{Z}^{d}$ with $0<|n_{*}|\leq N$, we have \begin{equation}\label{eq.reson}\|2\rho-\langle  n_{*},\tilde{\omega}\rangle\|_{\mathbb{R}/\mathbb{Z}}<\epsilon^{\frac{1}{15}},\end{equation}
then  $A_+(\mu)=e^{A''(\mu)}$ with $\lVert A''(\mu)\rVert \leqslant 2\epsilon^{\frac{1}{16}}.$ Moreover,
we have
$$
\rho(\tilde{\omega}, A_+(\mu)e^{ f_+(\mu,\cdot)} ) =\rho + \frac{ \langle  n_{*},\tilde{\omega} \rangle}{2}, \quad \text{for any} \quad  \mu \in \mathcal{O}.
$$

\end{enumerate}
\end{Lemma}

\begin{pf}
The proof follows \cite{CCYZ,LYZZ}. We will sketch the proof for completeness, but point out why we can make the conjugation continuous in the desired parameter set.

The basis of Lemma~\ref{eq.oneKam} is the following non-resonance cancellation lemma. We have a decomposition of the Banach algebra $\mathcal{B}_{\mathcal{O},r}(\mathrm{su}(1,1))$ into non-resonant spaces and
resonant spaces: $\mathcal{B}_{\mathcal{O},r}(\mathrm{su}(1,1)) = \mathcal{B}_{\mathcal{O},r}^{re}(\mathrm{su}(1,1),\eta) \oplus \mathcal{B}_{\mathcal{O},r}^{nre}(\mathrm{su}(1,1),\eta)$. Here, $\mathcal{B}_{\mathcal{O},r}^{nre}(\mathrm{su}(1,1),\eta)$ is defined in the following way: for any  $Y(\mu,\cdot)\in\mathcal{B}_{r}^{nre}(\mathrm{su}(1,1), \eta)$, we have
$$
A^{-1}(\mu)Y(\mu,\cdot+\tilde{\omega})A(\mu)\in \mathcal{B}_{\mathcal{O},r}^{nre}(\mathrm{su}(1,1), \eta)
$$
and
$$
\|A^{-1}(\mu)Y(\mu,\cdot+\tilde{\omega})A(\mu)-Y(\mu,\cdot)\|_{r}\geq\eta\|Y(\mu,\cdot)\|_{r}.
$$

\begin{Lemma}\label{HYLem}
Let $r>0$, $0 < \eta < 1$, $A(\mu)\in C^{0}(\mathcal{O},\mathrm{SU}(1,1))$, and $g(\mu,x) \in \mathcal{B}_{\mathcal{O},r}(\mathrm{su}(1,1))$. If
$$
\|g(\mu,\cdot)\|_{r}< \epsilon<   \frac{\eta^2}{C \sup_{\mu}\|A(\mu)\|^{4}},
$$
then there exist $Y(\mu,\cdot)\in\mathcal{B}_{\mathcal{O},r}(\mathrm{su}(1,1))$ and $g^{re}(\mu,\cdot)\in\mathcal{B}_{\mathcal{O},r}^{re}(\mathrm{su}(1,1),\eta)$ such that
$$
e^{-Y(\mu,x+\tilde{\omega})}A(\mu)e^{g(\mu,x)}e^{Y(\mu,x)} = A(\mu)e^{g^{re}(\mu,x)},
$$
with the estimate $\|Y(\mu,\cdot)\|_{r}\leq\epsilon^{\frac{1}{2}}$ and $\|g^{re}(\mu,\cdot)\|_{r}\leq 2\epsilon$.
\end{Lemma}

\begin{pf}
The proof is given in Appendix~\ref{AB} for the sake of curious readers.
\end{pf}

Once we have this key lemma, we can finish the proof. Denote the eigenvalues of a matrix $A$ by $\lambda(A)$. First note that we can always assume the eigenvalues of the constant matrix to be $e^{\pm i \rho}$. Otherwise, by continuity of the rotation number,
$$
\|\rho(\tilde{\omega}, A(\mu)e^{ f(\mu,\cdot)} ) -  \rho(\tilde{\omega}, A(\mu) )\|_0 \leq \epsilon,
$$
and we can modify $A(\mu)$, and write $\tilde{A}(\mu)e^{\tilde{f}(\mu,x)}=A(\mu)e^{f(\mu,x)}$, such that
\begin{equation}\label{cons-rho}
\lambda(\tilde{A}(\mu))= e^{\pm i\rho}
\end{equation}
and $\|\tilde{f}(\mu,\cdot)\|_{r}\leq  2 \epsilon$.
Then we can distinguish two cases:\\

\smallskip

\textbf{Non-resonant case:} In the non-resonant case (1), by noting \eqref{cons-rho}, similarly as in \cite[Proposition~3.1]{CCYZ}, one can  apply Lemma \ref{HYLem} once, and obtain $Y(\mu,\cdot)\in \mathcal{B}_{\mathcal{O},r}(\mathrm{su(1,1)})$ with $\|Y(\mu,\cdot)\|_{r}\leq 2\epsilon^{\frac{1}{2}}$ such that
$$
e^{-Y(\mu,x+\tilde{\omega})}\tilde{A}(\mu)e^{\tilde{f}(\mu,x)}e^{Y(\mu,x)}=\tilde{A}(\mu)e^{\tilde{f}^{re}(\mu,x)},
$$
where $\tilde{f}^{re}(\mu,\cdot)\in\mathcal{B}_{\mathcal{O},r}^{re}(\mathrm{su}(1,1),\epsilon^{\frac{1}{5}})$ with  $\| \tilde{f}^{re} \|_{r} \leq 4\epsilon$.

Furthermore, as in \cite[Proposition~3.1]{CCYZ}, one can easily check that for any $f \in \mathcal{B}_{\mathcal{O},r}^{re}(\mathrm{su}(1,1), \epsilon^{\frac{1}{5}})$, we have $\hat{f}(n)=0$ for all $0<|n|\leq N$. Therefore, by letting
$$
A_{+}=\tilde{A}(\mu)e^{\hat{\tilde{f}}^{re}(\mu,0)}, \qquad f_{+}=\sum_{|k|> N}\hat{\tilde{f}}^{re}(\mu,k)e^{i\langle k,x\rangle},
$$
where $\hat{\tilde{f}}^{re}(\mu,k)$ is the $k-$th Fourier coefficient of $\tilde{f}^{re}(\mu,\cdot)$, the proof in the non-resonant case is finished.\\

\smallskip

\textbf{Resonant case:} If \eqref{eq.reson} is satisfied, then by the fact $\tilde{\omega}\in DC(\tilde{\kappa},\tau)$, we have the following estimate,
$$
|2\rho| \geq \frac{\tilde{\kappa}}{|n_*|^{\tau}} -  \epsilon^{\frac{1}{15}} \geq \frac{\tilde{\kappa}}{ 2 |n_*|^{\tau}} \geq \frac{\tilde{\kappa}}{ 2 |N|^{\tau}}.
$$
On the other hand, by noting \eqref{cons-rho}, that is, the eigenvalues of $\tilde{A}(\mu)$ are fixed, we can apply Lemma~\ref{diagonalizing} to $\tilde{A}(\mu)$ and obtain $P(\mu)\in C^{0}( \mathcal{O}, \mathrm{SU}(1,1))$ with
$$
\|P(\mu)\|^{2}\leq \frac{2\|\tilde{A}(\mu)\|}{\rho}\leq \frac{1}{2}\epsilon^{-\frac{1}{150}}
$$
such that $$P^{-1}(\mu)\tilde{A}(\mu)P(\mu)=\begin{pmatrix}e^{i\rho} &0\\ 0&e^{-i\rho}\end{pmatrix} := A'(\mu).
$$

Let $g(\mu,\cdot)=P^{-1}(\mu)\tilde{f}(\mu,\cdot)P(\mu)$ with
$$
\|g(\mu,\cdot)\|_{r}\leq \|P(\mu)\|^{2}\|\tilde{f}(\mu,\cdot)\|_{r}\leq \frac{1}{4}\epsilon^{1-\frac{1}{75}}:=\epsilon'.
$$
Since $\epsilon^{\frac{1}{15}} \geq \|A'(\mu)\|^{2}\epsilon'$, we can apply Lemma~\ref{HYLem} with $\eta=\epsilon^{\frac{1}{15}}$ to the cocycle $(\tilde{\omega}, A'(\mu)e^{g(\mu,x)})$. We obtain $Y(\mu,\cdot)\in \mathcal{B}_{\mathcal{O},r}(\mathrm{su}(1,1))$ with $\|Y(\mu,\cdot)\|_{r}\leq\epsilon'^{\frac{1}{2}}$ such that
$$
e^{-Y(\mu,x+\tilde{\omega})}A'(\mu)e^{g(\mu,x)}e^{Y(\mu,x)}=A'(\mu)e^{g^{re}(\mu,x)}
$$
with $\|g^{re}(\mu,\cdot)\|_{r}\leq2\epsilon'$.

Define
\begin{equation}\label{eq.def1}
\mathcal{I}_{1}(\epsilon^{\frac{1}{15}})=\{n\in\mathbb{Z}^{d}:\|\langle n,\tilde{\omega}\rangle\|_{\mathbb{R}/\mathbb{Z}}\geq \epsilon^{\frac{1}{15}}\},
\end{equation}
\begin{equation}\label{eq.def2}
\mathcal{I}_{2}(\epsilon^{\frac{1}{15}})=\{n\in\mathbb{Z}^{d}:\|2\rho-\langle n,\tilde{\omega}\rangle\|_{\mathbb{R}/\mathbb{Z}}\geq\epsilon^{\frac{1}{15}}\}.
\end{equation}
Since any $F \in \mathcal{B}_{\mathcal{O},r}^{nre}(\mathrm{su}(1,1),\eta)$ takes the form
$$
\begin{aligned}&F(\mu,x)\\&=\sum_{n\in\mathcal{I}_{1}(\epsilon^{\frac{1}{15}})}\begin{pmatrix}
i\hat{t}(n) &0\\0 &-i\hat{t}(n)
\end{pmatrix}e^{i\langle n,x\rangle}+\sum_{n\in\mathcal{I}_{2}(\epsilon^{\frac{1}{15}})}\begin{pmatrix}
0 &\hat{v}(n)e^{i\langle n,x\rangle}\\
\overline{\hat{v}(n)}e^{-i\langle n,x\rangle} &0
\end{pmatrix},
\end{aligned}
$$
we may invoke the following results:
\begin{Claim}\cite{CCYZ,LYZZ}\label{res-set}
 Let $\mathcal{I}_{1}$ and $\mathcal{I}_{2}$ be given as in \eqref{eq.def1},\eqref{eq.def2}, then they have the following properties:
$$
\mathbb{Z}^{d}\backslash \mathcal{I}_{1}(\epsilon^{\frac{1}{15}})\bigcap\{n\in\mathbb{Z}^{d}: |n|\leq \kappa^{\frac{1}{\tau}}\epsilon^{-\frac{1}{15\tau}}\}=\{0\},
$$
$$
\mathbb{Z}^{d}\backslash\mathcal{I}_{2}(\epsilon^{\frac{1}{15}})\bigcap\{n\in\mathbb{Z}^{d}:|n|\leq 2^{-\frac{1}{\tau}}\kappa^{\frac{1}{\tau}}\epsilon^{-\frac{1}{15\tau}}-N\}=\{n_{*}\}.
$$
\end{Claim}

\medskip

Denote $N_{1}=2^{-\frac{1}{\tau}}\kappa^{\frac{1}{\tau}}\epsilon^{-\frac{1}{15\tau}}-N$. Then, as a consequence of Claim \ref{res-set}, we  have the following decomposition for $g^{re}(\mu,\cdot)$:
$$
\begin{aligned}g^{re}(\mu,\cdot)&=g_{0}^{re}(\mu,\cdot)+g_{n_{*}}^{re}(\mu,\cdot)+g_{2}^{re}(\mu,\cdot)\\
&=\begin{pmatrix}i\hat{t}(0) &0\\ 0&-i\hat{t}(0)\end{pmatrix}+\begin{pmatrix}0 &\hat{v}(n_{*})e^{i\langle n_{*},x\rangle}\\ \overline{\hat{v}(n)}e^{-i\langle n_{*},x\rangle} &0\end{pmatrix}\\
& \quad +\sum_{|n|\geq N_{1}}\hat{g}^{re}(\mu,n)e^{i\langle n,x\rangle}.
\end{aligned}
$$

Let $Q=\begin{pmatrix}e^{i\frac{\langle n_{*},x\rangle}{2}} &0\\ 0& e^{-i\frac{\langle n_{*},x\rangle}{2}}\end{pmatrix}$. Then
$$
Q^{-1}(x+\tilde{\omega})A'(\mu)e^{g^{re}(\mu,x)}Q(x)=\tilde{A}'(\mu)e^{\tilde{g}(\mu,x)},
$$
where
$$
\tilde{A}'(\mu)=\begin{pmatrix}e^{i(\rho-\langle n_{*},\tilde{\omega}\rangle/2)} &0\\0 &e^{-i(\rho-\langle n_{*},\tilde{\omega}\rangle/2)}\end{pmatrix}
$$
and $\tilde{g}(\mu,\cdot)=Q^{-1}g^{re}(\mu,\cdot)Q$. Then let
$$
A_{+}=\tilde{A}'(\mu)\exp\begin{pmatrix}i\hat{t}(0) &\hat{v}(n_{*})e^{i\langle n_{*},x\rangle}\\ \overline{\hat{v}(n)}e^{-i\langle n_{*},x\rangle} &-i\hat{t}(0)\end{pmatrix}
$$
and $F(\mu,\cdot)=Q^{-1}g_{2}^{re}(\mu,\cdot)Q$; $f_{+}$ is given by
$$
e^{g^{re}_{0}(\mu,\cdot)+g^{re}_{n_{*}}(\mu,\cdot)+F(\mu,\cdot)}=e^{g^{re}_{0}(\mu,\cdot)+g^{re}_{n_{*}}(\mu,\cdot)}e^{f_{+}}
$$
with the estimate $\|f_{+}\|_{r'}\leq\epsilon^{2}$. Moreover, the conjugation given by the above argument is $Z(\mu,\cdot)= P(\mu) \cdot e^{Y(\mu,\cdot)} \cdot Q\in \mathcal{B}_{\mathcal{O},r'}(\mathrm{SU}(1,1))$. We should point out that the index $n_{*}$ in the rotation $Q$ is independent of the choice of $\mu\in\mathcal{O}$, which is key to our proof.
\end{pf}

We prove Proposition~\ref{conDep} by iteration of Lemma~\ref{eq.oneKam}. In the application, we just take $\tilde{\omega}=2\omega$ and $\mathcal{O}=\mathcal{R}_{\tilde{\rho}}.$ Without loss of generality, we assume that $r<1$ and $\tilde{r}=\frac{r}{2}$.  Note that since we take $\mu\in[0, \frac{1}{2}]\times\partial\mathbb{D}$, we
can take $|\delta|$ small enough such that
\begin{eqnarray*}
\|T_{\delta}(\mu,\cdot)-T_{0}(\mu)\|_{r}  & \leq&  \frac{1}{\sqrt{(1-\lambda_{1}^{2})(1-\lambda_{2}^{2})}} |\delta| \|h\|_{r} \\ &\leq&  \epsilon_* \leq        \frac{D_0}{\sup_{\mu}\|T_{0}(\mu)\|^{C_{0}}}(\frac{r-\tilde{r}}{8})^{C_0\tau},
\end{eqnarray*}
where $D_0=D_0(\tilde{\kappa},\tau,d)$ is the constant defined in Lemma~\ref{eq.oneKam}. By the Implicit Function Theorem, we can always write
$$
T_{\delta}(\mu,\cdot) = T_{0}(\mu) e^{f(\mu,\cdot)} :=  A(\mu)e^{f(\mu,\cdot)}
$$
with the estimate $\|f(\mu,\cdot)\|_{r} \leq \epsilon_*$. Then we can define the sequence inductively:
$$
\epsilon_{j} = \epsilon_*^{2^{j}}, \quad r_j-r_{j+1}=\frac{r-\frac{r+\tilde{r}}{2}}{4^{j+1}}, \qquad N_j=\frac{2|\ln\epsilon_j|}{r_j-r_{j+1}}.
$$
Assume that we are at the $(j+1)^{th}$ KAM step, that is, we have already constructed $Z_j \in C^\omega_{r_{j}}( \mathcal{R}_{\tilde{\rho}} \times 2 \T^d, \mathrm{SU}(1,1))$ such that
$$
Z^{j}(\mu,x+2\omega)^{-1}A(\mu)e^{f(\mu,x)}Z^{j}(\mu,x)=A_{j}(\mu)e^{f_{j}(\mu, x)}
$$
with $\|f_j\|_{r_j}\leq \epsilon_j$ and
\begin{equation}\label{rhoj}
2\rho_{j} = \rho(2\omega, A_{j}(\mu)e^{f_{j}(\mu,x)} )= \tilde{\rho} +\sum_{i=1}^{j}   \langle    n_{i}, 2\omega \rangle  \mod \Z   \;   \text{ for any } \mu \in \mathcal{R}_{\tilde{\rho}}.
\end{equation}

By our selection of  $\epsilon_0$, one can check that
\begin{equation}\label{iter}
\epsilon_j \leq \frac{D_0}{\sup_{\mu}\|A(\mu)\|^{C_{0}}}(r_j-r_{j+1})^{C_0\tau}.
\end{equation}
Indeed, $\epsilon_j$ on the left side of the inequality decays at least super-exponentially in $j$, while $(r_j-r_{j+1})^{C_0\tau}$ on the right side decays exponentially in $j$.

Note that \eqref{iter} implies that Proposition~\ref{eq.oneKam} can be applied iteratively. Consequently, one can construct
$$
\tilde{Z}^{j+1}(\mu,x+2\omega)^{-1}A_{j}(\mu)e^{f_{j}(\mu,x)}\tilde{Z}^{j+1}(\mu,x)=A_{j+1}(\mu)e^{f_{j+1}(\mu,x)}
$$
with $\|f_{j+1}\|_{r_{j+1}}\leq \epsilon_{j+1}$. More precisely, we can distinguish two cases:\\

\smallskip
\textbf{Non-resonant case:} If for any $n\in\mathbb{Z}^{d}$ with $0<|n|\leq N_j$, we have
$$
\|2\rho_j-\langle n,2\omega\rangle\|_{\mathbb{R}/\mathbb{Z}}\geq\epsilon_j^{\frac{1}{15}},
$$
then we have $\| \tilde{Z}^{j+1}(\mu,\cdot)-Id\|_{r_{j+1}}\leq \epsilon_j^{\frac{1}{2}}$ and
$$
\rho(2\omega, A_{j+1}(\mu)e^{f_{j+1}(\mu,x)} ) =\rho_j, \quad \text{for any} \quad  \mu \in \mathcal{R}_{\tilde{\rho}}.
$$
In this case, we just set $n_{j+1}=0$. \\

\smallskip
\textbf{Resonant case:} If for some $n_{*}\in\mathbb{Z}^{d}$ with $0<|n_{*}|\leq N_j$, we have
\begin{equation*}
\|2\rho_j-\langle  n_{*},2\omega\rangle\|_{\mathbb{R}/\mathbb{Z}}<\epsilon_j^{\frac{1}{15}},
\end{equation*}
then we have
$$
2\rho(2\omega, A_{j+1}(\mu)e^{f_{j+1}(\mu,x)} )  =2\rho_j + \langle  n_{*},2\omega \rangle \mod \Z \; \text{ for any } \mu \in \mathcal{R}_{\tilde{\rho}}.
$$
In this case, we just set $n_{j+1}=n_{*}$.

\medskip

By letting $Z^{j+1}(\mu,\cdot)=Z^{j}(\mu,\cdot)\cdot\tilde{Z}^{j+1}(\mu,\cdot)$, we obtain
$$
Z^{j+1}(\mu,x+2\omega)^{-1}A(\mu)e^{f(\mu,x)}Z^{j+1}(\mu,x)=A_{j+1}e^{f_{j+1}(x)},
$$
and  by \eqref{rhoj}, we always have
$$
2\rho_{j+1} := \rho(2\omega, A_{j+1}(\mu)e^{f_{j+1}(\mu,x)} )= \tilde{\rho} +\sum_{i=1}^{j+1}   \langle    n_{i}, 2\omega \rangle  \mod \Z \; \text{ for any } \mu \in \mathcal{R}_{\tilde{\rho}}.
$$

According to \cite[Lemma~3]{Eli92}, when the rotation number is Diophantine or rational with respect to the frequency, resonance happens only finitely many times, that is, if $j$ is large enough, we always have  $\| \tilde{Z}^{j+1}(\mu,\cdot)-Id\|_{r_{j+1}}\leq \epsilon_j^{\frac{1}{2}}$, and $n_{j+1}=0$. Then $Z(\mu,x)=\lim\limits_{j\to\infty} Z^{j}(\mu,x)$ serves our purpose.

\subsection{Cantor Spectrum}

We also need to establish a result similar to Proposition~\ref{openGap} for the model considered in this section. Note that for the same reason as above, we can focus on the case $k \not= 0$.

\begin{Proposition}\label{openGapLin2}
Denote the tongue boundaries of $(2\omega, T_{\delta}(\lambda_2,z))$  with label $k\neq 0$ by $\theta_{k}^{\pm}(\delta)$. If
$|\lambda_{2}| < r_{2}(\lambda_{1})$ is small enough, we have
$$
\frac{d\theta_{k}^{+}}{d\delta}(0)\neq \frac{d\theta_{k}^{-}}{d\delta}(0)$$
as long as $\hat{h}_{k}\neq 0$.
\end{Proposition}

\begin{pf}
Let $\delta_{0}=0$ and $z_{0}=e^{i\theta_{0}}$. Then the matrix $T_{0}(\lambda_2,z)$ has the form
$$
T_{0}(\lambda_2,e^{i\theta_{0}})=\frac{1}{\sqrt{(1-\lambda_{1}^{2})(1-\lambda_{2}^{2})}}\begin{pmatrix}e^{i\theta_{0}}+\lambda_{1}\lambda_{2} &-(\lambda_{1}e^{-i\theta_{0}}+\lambda_{2})\\-(\lambda_{1}e^{i\theta_{0}}+\lambda_{2}) &e^{-i\theta_{0}}+\lambda_{1}\lambda_{2}\end{pmatrix}.
$$
As we are considering tongue boundaries of the gap with label $k$, the eigenvalues of $T_0$ are given by $\lambda_{\pm}(\theta) = e^{\pm i\langle k,\omega\rangle}$. Indeed, by direct computation, one can show that the eigenvalues are given by
\begin{equation}\label{eigen}
\lambda_{\pm}=\frac{\lambda_{1}\lambda_{2}+\cos\theta_{0}}{\sqrt{(1-\lambda_{1}^{2})(1-\lambda_{2}^{2})}}\pm i\sqrt{1-\frac{(\lambda_{1}\lambda_{2}+\cos\theta_{0})^{2}}{(1-\lambda_{1}^{2})(1-\lambda_{2}^{2})}}.
\end{equation}
Here we should point that for any given $k$, there exist two values of $\theta_0$ that share the same eigenvalue, this is because $(2\omega, T_{\delta}(\lambda_2,e^{i\theta_{0}}))$ is the twice iterate of the cocycle $(T_{\omega}, S(\alpha,e^{i\theta_{0}}))$.

Let $U=\frac{1}{\sqrt{\cos2\tilde{\theta}}}\begin{pmatrix}\cos\tilde{\theta} &e^{2i\phi}\sin\tilde{\theta}\\ e^{-2i\phi}\sin\tilde{\theta} &\cos\tilde{\theta}\end{pmatrix}$ be given by Lemma~\ref{diagonalizing} such that
$$
U^{-1}T_{0}U=\begin{pmatrix}e^{i\langle k,\omega\rangle} &0\\0 &e^{-i\langle k,\omega\rangle}\end{pmatrix}.
$$
Let
$$
Z(x)=U\begin{pmatrix}e^{i\frac{\langle k,x\rangle}{2} }&0\\ 0&e^{-i\frac{\langle k,x\rangle}{2}}\end{pmatrix} =\frac{1}{\sqrt{\cos2\tilde{\theta}}}\begin{pmatrix}
\cos\tilde{\theta}e^{i\frac{\langle k,x\rangle}{2}} &e^{2i\phi}\sin\tilde{\theta}e^{-i\frac{\langle k,x\rangle}{2}}\\
e^{-2i\phi}\sin\tilde{\theta}e^{i\frac{\langle k,x\rangle}{2}} &\cos\tilde{\theta}e^{-i\frac{\langle k,x\rangle}{2}}
\end{pmatrix}.
$$
Then we have
$$
Z^{-1}(x+2\omega)T_{0}(\lambda_2,e^{i\theta_{0}})Z(x)=Id,
$$
and it conjugates the perturbed cocycle $(2\omega,T_{\delta}(\lambda_2,e^{i(\theta_{0}+\eta)}))$ into
$$
\begin{aligned}
&Z^{-1}(x+2\omega)T_{\delta}(\lambda_2,e^{i(\theta_{0}+\eta)})Z(x)\\
&=Z^{-1}(x)T^{-1}_{0}(\lambda_2,e^{i\theta_{0}})T_{\delta}(\lambda_2,e^{i(\theta_{0}+\eta)})Z(x)=e^{\tilde{P}(x,\zeta)+\tilde{O}_{2}(x,\zeta)}.
\end{aligned}
$$
According to \eqref{Peq1} with $\delta_{0}=0$, we have
\begin{eqnarray*}
\tilde{P} =  \begin{pmatrix} i a_1 (\lambda_2,x)\eta+ i b_1(\lambda_2,x)\beta&a_2 (\lambda_2,x)\eta+b_2(\lambda_2,x)\beta \\\overline{a_2(\lambda_2,x)}\eta+\overline{b_2(\lambda_2,x)}\beta& - ia_1(\lambda_2,x)\eta-i b_1( \lambda_2,x)\beta \end{pmatrix} ,
\end{eqnarray*}
where
\begin{eqnarray*}
a_2(\lambda_2,x) &= &\frac{i\left(e^{-i\langle k,x\rangle}(e^{2i\phi}\sin2\tilde{\theta}-\lambda_{2}e^{i(\theta_{0}+4\phi)}\sin^{2}\tilde{\theta}-\lambda_{2}\cos^{2}\tilde{\theta}e^{-i\theta_{0}}\right)}{(1-\lambda_{2}^{2})\cos2\tilde{\theta}},\\
b_2(\lambda_2,x) &=& \frac{ie^{-i\langle k,x\rangle}(2e^{2i\phi}\sin\tilde{\theta}\cos\tilde{\theta}f_{11}+e^{4i\phi}\sin^{2}\tilde{\theta}\overline{f_{12}}+\cos^{2}\tilde{\theta}f_{12})}{(1-\lambda_{2}^{2})\cos2\tilde{\theta}},
\end{eqnarray*}
where $f_{11}$, $f_{12}$ are  defined in \eqref{eq.coefficients3}--\eqref{eq.coefficients4}. In this case, they take the form
\begin{eqnarray*}
f_{11} (\lambda_2)&=&\lambda_{1}^{2}\lambda_{2}^{2}(h_{+}-h)-\lambda_{1}^{2}h-\lambda_{2}^{2}h_{+}-2\lambda_{1}\lambda_{2}h\cos \theta_{0},
\end{eqnarray*}
\begin{eqnarray*}
f_{12} (\lambda_2)&=&(\lambda_{1}e^{-2i\theta_{0}}+2\lambda_{1}^{2}\lambda_{2}e^{-i(\theta_{0})} +\lambda_1\lambda_2^2 ) h
   +  \lambda_{2} (1-\lambda_1^2)e^{-i\theta_{0}}  h_{+}.
\end{eqnarray*}

By taking averages, we obtain $a_2(\lambda_2) = [a_2(\lambda_2,x)]=0$ and
$$
b_2(\lambda_2) =  \frac{i}{(1-\lambda_{2}^{2})\cos2\tilde{\theta}}   \left( (2e^{2i\phi}\sin\tilde{\theta}\cos\tilde{\theta}\tilde{f}_{11}+e^{4i\phi}\sin^{2}\tilde{\theta}\overline{\tilde{f}_{12}}+\cos^{2}\tilde{\theta}\tilde{f}_{12})\hat{h}_{-k}\right),
$$
where we simply denote $[e^{-i\langle k,x\rangle}f_{ij}]=\tilde{f}_{ij}\hat{h}_{-k}$ for $j=1,2$. By Corollary~\ref{openLinearly} (2),(3), there exist two analytic tongue boundaries $\theta^{\pm}_{k}(\delta)$, and $\frac{d\theta_{k}^{+}}{d\delta}(0)\neq\frac{d\theta_{k}^{-}}{d\delta}(0)$ if and only if $b_{2}(\lambda_2)\neq 0$.

Just note
\begin{eqnarray*}
 \overline{b_2(\lambda_2)}  &=& \frac{ -i\hat{h}_{k}(e^{-4i\phi}\sin^{2}\tilde{\theta}\tilde{f}_{12}+2e^{-2i\phi}\sin\tilde{\theta}\cos\tilde{\theta} \tilde{f}_{11}+\cos^{2}\tilde{\theta}\overline{\tilde{f}_{12}})}{(1-\lambda_{2}^{2})\cos2\tilde{\theta}} \\
   &=&  \frac{-i\hat{h}_{k}\cos^{2}\tilde{\theta}\overline{\tilde{f}_{12}}\left(\left(e^{-2i\phi}\tan\tilde{\theta}\frac{\tilde{f}_{12}}{|\tilde{f}_{12}|}+\frac{\tilde{f}_{11}}{|\tilde{f}_{12}|}\right)^{2}+1-\left(\frac{ \tilde{f}_{11}}{|\tilde{f}_{12}|}\right)^{2}\right) } {(1-\lambda_{2}^{2})\cos2\tilde{\theta}}.
\end{eqnarray*}

A key observation here is that there exists a positive number $r_{2}(\lambda_{1})$ such that if $|\lambda_{2}|<r_{2}(\lambda_{1})$, we have
\begin{equation}\label{eq.smallCond}
\left|\frac{\tilde{f}_{11}(\lambda_2)}{\tilde{f}_{12}(\lambda_2)}\right|<1
\end{equation}
To see this, by taking $\lambda_{2}=0$, we have $\tilde{f}_{11}(0)=-\lambda_{1}^{2}$ and $\tilde{f}_{12}(0)=\lambda_{1}e^{-2i\theta_{0}}$,
which implies that \eqref{eq.smallCond} holds for $\lambda_2=0$. Then  the conclusion follows via a simple continuity argument.

According to \eqref{eq.smallCond} and our choice of $\lambda_{2}$, and note according to Lemma \ref{diagonalizing}, $\tilde{\theta}\in(-\frac{\pi}{2},\frac{\pi}{2})$ never takes the value $\frac{\pi}{4}$, $b_2(\lambda_2)$ does not vanish if $\hat{h}_{k}\neq 0$. This completes the proof.
\end{pf}

\medskip

\noindent\textbf{Proof of  Theorem \ref{MainThm1}:}
For $k\in\mathbb{Z}^{d}$, let $\mathcal{G}_{k}$ be the collection of $h\in C^{\omega}(\T^{d},\mathbb{R})$ such that the $k$-th Fourier coefficient of $h$ satisfies $\hat{h}_{k}\neq 0$. Define $\mathcal{G}=\cap_{k\in\mathbb{Z}^{d}}\mathcal{G}_{k}$, clearly, $\mathcal{G}$ is generic. By Theorem \ref{Thm.period2}, the tongue boundaries $\theta_k^{\pm}(\delta)$ of $(2\omega,T_{\delta}(\lambda_2,e^{i\theta}))$ are analytic function of $\delta$. By Proposition \ref{Lem.gapLabel}, there exists two  branches  $\theta^{\pm}_{k,i}(\delta)$ $i=1,2$ of the tongue boundary
$(T_{\omega},S(\alpha,z))$, according to its rotation number   \begin{equation*}
2\rho(T_\omega,S(\alpha,e^{i\theta}))=\langle k,\omega\rangle \mod \mathbb{Z},
\end{equation*}
or
\begin{equation*}
2\rho(T_\omega,S(\alpha,e^{i\theta}))=\langle k,\omega\rangle +\frac{1}{2} \mod \mathbb{Z}.
\end{equation*}
Clearly, they are both real analytic. Indeed, in the case $\delta=0$, this can be clearly shown in \eqref{eigen}. As a consequence of Proposition \ref{openGapLin2},  these tongue boundaries have a transversality at the origin $\delta=0$, thus the set of couplings $|\delta|<\epsilon$ such that certain tongue boundaries coincide is finite. We thus finish the proof.
\qed

\medskip

\noindent\textbf{Proof of  Theorem \ref{Thm.QW}:}

In the case $C_{n}\in \mathrm{SU}(1,1)$, which means that $c_{n}^{11}=\overline{c}_{n}^{22}$ and $c_{n}^{12}=\overline{c}_{n}^{21}$. In this scenario, we have $\sigma_{n}^{1}=-\sigma_{n}^{2}, n\in \mathbb{Z}$ and therefore the Verblunsky coefficients given by \eqref{eq.VC} take the form
$$
\alpha_{2n+1}=0,\alpha_{2n+2}=c_{n+1}^{12}, n\in \mathbb{Z}.
$$
In order to make Theorem~\ref{MainThm1} applicable in the present setting, let $c_{n}^{12}=\lambda e^{i h(x+n\omega)}$, where $h,\omega$ are as in Theorem \ref{MainThm1} and $\lambda\in (0,1)$.
Note that our discussion in Section~\ref{AlmostModel} shows that $\lambda_{1},\lambda_{2}$ are in a symmetric position. Thus by taking $\lambda_{1}=\lambda, \lambda_{2}=0$  in \eqref{eq.sampling}--\eqref{eq.way}, we recover $\alpha_{2n+1}=0,\alpha_{2n+2}=\lambda e^{ih(x+(n+1)2\omega)}$. One can easily align this extra factor $2$ by adjusting either side. Then   Theorem \ref{Thm.QW} follows immediately from Theorem \ref{MainThm1}.  \qed

\appendix

\section{Proof of Theorem \ref{KAM}}\label{appendix1}

In this section, we give a proof of Theorem~\ref{KAM}, the proof of Theorem \ref{KAM1} is very similar. This proof follows the same lines as the proofs of Theorems 4 and 11 in \cite{Puig06}, although they are given in different groups.

Let $A_{0}=\chi^{p}\begin{pmatrix}0 &1\\ 0 &0\end{pmatrix}$ or $0$ with $\chi\neq 0$. For any fixed perturbation $\chi^{p+1}P(x,\mu)$, we are going to find a counterterm $\chi^{p+1}\xi^{*}(\mu)$ and a conjugation $Z(x)=\exp(\chi X(x))$ such that
\begin{equation}
Z^{-1}(x+\omega)e^{A_{0}}e^{\chi^{p+1}P(x,\mu)}e^{-\chi^{p+1}\xi^{*}(\mu)}Z(x)=e^{A_{0}}.
\end{equation}
We omit the factor $\chi^{p+1}$ since it does not affect the iteration process and still use the notations $P_{j},\eta_{j}$ for $\chi^{p+1}P_{j},\chi^{p+1}\eta_{j}$.
In order to do this, we construct the following iterative process:
\begin{equation}\label{cohomoeq4}
Z^{j}(x+\omega)^{-1}e^{A_{0}}e^{P_{j}(x,\mu)}e^{-\eta_{j}(\mu)}Z^{j}(x)=e^{A_{0}}e^{P_{j+1}(x)}e^{-\eta_{j+1}(\mu)},
\end{equation}
where $P_{j+1},\eta_{j+1}$ is smaller than $P_{j},\eta_{j}$. The initial choice of $P_{0}$ is $P$ and $\eta_{0}$ is to be defined along the iterative process. The relation between $\eta_{j}$ and $\eta_{j+1}$ is given by $\xi^{j}:\eta_{j+1}\to\eta_{j}$ such that $\eta_{j+1}=\xi^{j}(\eta_{j+1})-C([P_{j}(\xi^{j}(\eta_{j+1}))])$, or equivalently $\eta_{j+1}=\eta_{j}-C([P_{j}(\eta_{j})])$. Define $Z_{j}=Z^{j}\circ Z_{j-1}$ and $\xi_{j}=\xi_{j-1}\circ\xi^{j}$.
If $Z_{j},\xi_{j}$ converge to some $Z,\xi^{*}$, then we have $$Z^{-1}(x+\omega)e^{A_{0}}e^{\chi^{p}P(x,\mu)}e^{-\chi^{p}\xi^{*}(0)}Z(x)=e^{A_{0}}.$$
Now consider the linearized equation of \eqref{cohomoeq4}, replace all the $e^{X}$ terms by $Id+X$ and denote the error terms $e^{X}-Id-X$ by $E(X)$.
Then we have the following expressions for the perturbation of the next step:
\begin{equation}\label{errorterm}
\begin{aligned}&P_{j+1}-P_{j+1}\eta_{j+1}+E(P_{j+1})e^{-\eta_{j+1}}+e^{P_{j+1}}E(-\eta_{j+1})\\
&=A_{0}^{-1}E(-X_{j}(x+\omega))A_{0}e^{P_{j}}e^{-\eta_{j}}Z^{j}+e^{P_{j}}e^{-\eta_{j}}E(X_{j})\\
&+(P_{j}-\eta_{j})X_{j}(x)-P_{j}\eta_{j}X_{j}+E(P_{j})e^{\eta_{j}}Z^{j}+e^{P_{j}}E(-\eta_{j})Z^{j}\\
&-A_{0}^{-1}X_{j}(x+\omega)A_{0}(X_{j}+P_{j}-\eta_{j})-A_{0}^{-1}X_{j}(x+\omega)A_{0}(e^{P_{j}}e^{-\eta_{j}}E(X_{j})\\
&+(P_{j}-\eta_{j})X_{j}-P_{j}\eta_{j}X_{j}+E(P_{j})e^{-\eta_{j}}Z_{j}+e^{P_{j}}E(-\eta_{j})Z^{j})
\end{aligned}
\end{equation}

\begin{Lemma}(The inductive lemma)
Assume $(A_{0},C,S,\omega)$ is admissible with constants $c'$ and $\nu$. Fix a complex domain $$D_{j}:\quad |\im x|<\rho_{j},\quad |\eta_{j}|<\sigma_{j}$$
and a constant $0<\delta_{j}<\rho_{j}$. Then there exists a constant $K$ such that if $P_{j}$ is analytic on $D_{j}$ and belongs to $\sl(2,\mathbb{R})$ for real values $(x,\eta_{j})$, and
$$
\|P_{j}\|_{D_{j}}= \sup_{(x,\eta_{j})\in D_{j}}\|P_{j}(x,\eta_{j})\|\leq \epsilon_{j}\leq K\sigma_{j},
$$
then in the domain
$$
D_{j+1}:\quad \|\im x\|< \rho_{j}-\delta_{j},\quad \|\eta_{j}\|<\frac{\sigma_{j}}{2}
$$
the transformation
$$
Z^{j}(x)=\exp( \chi X_{j}(x)),
$$
where $X_{j}(x)$ satisfies
$$
-A_{0}^{-1} X_{j}(x+\omega)A_{0}+X_{j}(x)+P_{j}(x,\eta_{j})-C([P_{j}(x,\eta_{j})])=0,[X_{j}]=S([P_{j}]),
$$
is real analytic, and the equation
\begin{equation}\label{eta}
\eta_{j+1}=\eta_{j}-C([P_{j}(x,\eta_{j})])
\end{equation}
defines a local diffeomorphism $\xi^{j}:  \eta_{j+1}\in D(0,\epsilon_{j})\to \xi^{j}(\eta_{j+1})\in D(0,2\epsilon_{j})$ such that the equation
$$
Z^{j,-1}(x+\omega)A_{0}e^{P_{j}(x,\xi^{j}(\eta_{j+1}))}e^{-\xi^{j}(\eta_{j+1})}Z^{j}(x)=A_{0}e^{P_{j+1}(x,\eta_{j+1})}e^{-\eta_{j+1}}
$$
holds in the domain
$$
D^{j+1}: \quad |\im x|<\rho_{j}-\delta_{j},\quad |\eta_{j+1}|<\epsilon_{j}
$$
with the following estimates:
\begin{equation}\label{estX}
\|X_{j}\|_{D_{j+1}}\leq M:=c\frac{\epsilon_{j}}{\delta_{j}^{\nu}},
\end{equation}
\begin{equation}\label{estP}
\|P_{j+1}\|_{D^{j+1}}\leq \epsilon_{j+1}:=\epsilon_{j}^{\frac{3}{2}},
\end{equation}
\begin{equation}\label{estDeta}
\|D_{\eta_{j+1}}\xi^{j}\|_{\epsilon_{j}}\leq 1+c_{1}\frac{\epsilon_{j}}{\sigma_{j}}.
\end{equation}
\end{Lemma}

\begin{pf}
\textit{Bounds for $P_{j+1}$}: Recall that $P_{j+1}$ can be expressed as in \eqref{errorterm}. In order to estimate $\|P_{j+1}\|_{D^{j+1}}$, we seek to estimate \eqref{errorterm} term by term. Note that
$$
\|P_{j}\|_{D_{j}}\leq \epsilon_{j},\|X_{j}\|_{D_{j}}\leq c\frac{\epsilon_{j}}{\delta_{j}^{\nu}},\|\eta_{j+1}\|\leq \epsilon_{j},\|\eta_{j}\|\leq 2\epsilon_{j}.
$$
Thus, $\|E(X_{j})\|\sim \|X_{j}\|^{2}$ gives
$$
\|A_{0}^{-1}E(-X_{j}(x+\omega))A_{0}e^{P_{j}}e^{-\eta_{j}}Z^{j}\|_{D^{j+1}}\sim \|X_{j}\|^{2}_{D_{j}}\leq c^{2}\frac{\epsilon_{j}^{2}}{\delta_{j}^{2\nu}}
$$
and
$$
\|e^{P_{j}}e^{-\eta_{j}}E(X_{j})\|_{D^{j+1}}\sim\|X_{j}\|^{2}_{D_{j}}\leq c^{2}\frac{\epsilon_{j}^{2}}{\delta_{j}^{2\nu}},
$$
$$
\|(P_{j}-\eta_{j})X_{j}(x)-P_{j}\eta_{j}X_{j}\|_{D^{j+1}}\leq \frac{\epsilon^{2}_{j}}{\delta^{\nu}_{j}}.
$$
Moreover, we have
$$
\|e^{P_{j+1}}E(-\eta_{j+1})\|_{D^{j+1}}\leq 2\|\eta_{j+1}\|^{2}\leq 2\epsilon^{2}_{j},
$$
$$
\|E(P_{j})e^{\eta_{j}}Z^{j}+e^{P_{j}}E(-\eta_{j})Z^{j}\|_{D^{j+1}}\sim\|P_{j}\|^{2}_{D_{j}}+\|\eta_{j}\|^{2}_{D_{j}}\leq 5\epsilon_{j}^{2},
$$
$$
\|A_{0}^{-1}X_{j}(x+\omega)A_{0}(X_{j}+P_{j}-\eta_{j})\|_{D^{j+1}}\leq 3\|X_{j}\|_{D_{j}}^{2}\leq 3c^{2}\frac{\epsilon^{2}_{j}}{\delta_{j}^{2\nu}},
$$
$$
\|A_{0}^{-1}X_{j}(x+\omega)A_{0}(e^{P_{j}}e^{-\eta_{j}}E(X_{j})+(P_{j}-\eta_{j})X_{j})\|_{D^{j+1}}\leq 3c^{3}\frac{\epsilon^{3}_{j}}{\delta_{j}^{3\nu}},
$$
$$
\|P_{j}\eta_{j}X_{j}\|_{D^{j+1}}\leq 2c\frac{\epsilon^{3}_{j}}{\delta^{\nu}_{j}},
$$
$$
\|A_{0}^{-1}X_{j}(x+\omega)A_{0}(E(P_{j})e^{-\eta_{j}}Z_{j}+e^{P_{j}}E(-\eta_{j})Z^{j})\|_{D^{j+1}}\leq 5c\frac{\epsilon^{3}_{j}}{\delta_{j}^{\nu}}.
$$
We put these together and obtain
$$
\|P_{j+1}\|_{D^{j+1}}\leq \epsilon^{2}_{j}(\frac{5c^{2}}{\delta_{j}^{2\nu}}+\frac{1}{\delta_{j}^{\nu}}+7+\frac{3c^{3}\epsilon_{j}}{\delta_{j}^{3\nu}}+\frac{7c\epsilon_{j}}{\delta_{j}^{\nu}}).
$$
In order to get $\|P_{j+1}\|_{D^{j+1}}\leq \epsilon_{j}^{\frac{3}{2}}$, we need to pick $\epsilon_{0}$ such that
\begin{equation}\label{cond1}
\epsilon^{\frac{1}{2}}_{j} \left(\frac{5c^{2}}{\delta_{j}^{2\nu}}+\frac{1}{\delta_{j}^{\nu}}+7+\frac{3c^{3}\epsilon_{j}}{\delta_{j}^{3\nu}}+\frac{7c\epsilon_{j}}{\delta_{j}^{\nu}} \right)\leq 1.
\end{equation}
\textit{Inversion of \eqref{eta}}: Let $F^{j}(\eta_{j})=C([P_{j}(\eta_{j})])$. Then $F^{j}$ is analytic on the ball $D(0,\sigma_{j})$. By the Cauchy inequality we have
$$
\|D_{\eta_{j}}F^{j}(\eta_{j})\|_{\frac{\sigma_{j}}{2}}\leq c'\frac{\|F^{j}\|_{\sigma_{j}}}{\sigma_{j}-\frac{\sigma_{j}}{2}}\leq 2c'\frac{\epsilon_{j}}{\sigma_{j}},
$$
where $c'$ is a universal constant. Assume that $K\leq \min \{\frac{1}{4},\frac{1}{4c'}\}$. We want to show that the map $\xi^{j}$ given by
$$
\eta_{j+1}\in D(0,\epsilon_{j})\to \xi^{j}(\eta_{j+1})\in D \left(0,\frac{\sigma_{j}}{2} \right)
$$
is well defined. Since $\epsilon_{j}\leq K\sigma_{j}$, we have $2\epsilon_{j}\leq \frac{\sigma_{j}}{2}$, and therefore
$$
\|\xi^{j}(\eta_{j+1})\|\leq \|\eta_{j+1}\|+\|F^{j}(\eta_{j})\|\leq 2\epsilon_{j}\leq \frac{\sigma_{j}}{2}.
$$
Moreover, we have
$$
\|D_{\eta_{j+1}}\xi^{j}\|_{\epsilon_{j}}\leq \frac{1}{1-\|D_{\eta_{j}}F^{j}\|_{\frac{\sigma_{j}}{2}}} \leq\frac{1}{1-2c'\frac{\epsilon_{j}}{\sigma_{j}}}\leq 1+c_{1}\frac{\epsilon_{j}}{\sigma_{j}},
$$
where $c_{1}=4c'$. This verifies \eqref{estDeta}.

In order to perform this iterative process to any order, we need to take $$
\epsilon_{j}=\epsilon_{0}^{(\frac{3}{2})^j}, \, \sigma_{j+1}=\epsilon_{j}, \, \rho_{j}=\rho_{0} \left(\frac{1}{2}+\frac{1}{2^{j+1}} \right), \, \rho_{j+1}=\rho_{j}-\delta_{j}, \, \delta_{j}=\frac{\rho_{0}}{2^{j+1}}.
$$
\textit{Choice of $\epsilon_{0}$:} In the one step process, there are three conditions for the choice of $\epsilon_{0}$, namely:
\begin{equation}\label{con1}
\epsilon_{j}<K\sigma_{j},
\end{equation}
\begin{equation}\label{con2}
\|X_{j}\|_{D_{j+1}}\leq c\frac{\epsilon_{j}}{\delta_{j}^{\nu}}<\epsilon_{j}^{\frac{1}{2}},
\end{equation}
and \eqref{cond1}. Condition \eqref{con1} is equivalent to $\epsilon_{0}^{\frac{1}{2}(\frac{3}{2})^{j-1}}<K,$
which is true for all $j\geq 0$ if we pick
\begin{equation}\label{con11}
\epsilon_{0}<K^{3}.
\end{equation}
Condition \eqref{con2} is equivalent to $\epsilon_{0}^{\frac{1}{2}(\frac{3}{2})^{j}}<\frac{1}{c}(\frac{\rho_{0}}{4})^{\nu}\frac{1}{2^{j\nu}},$ which is true if we choose $\epsilon_{0}$ with
\begin{equation}\label{con22}
\epsilon_{0} < \min \left\{\frac{1}{c^{2}} \left(\frac{\rho_{0}}{4} \right)^{2\nu},\exp \left(2\nu\log\frac{1}{4} \right) \right\}.
\end{equation}
In order to have \eqref{cond1}, we need to pick $\epsilon_{0}$ with
\begin{equation}\label{cond11}
\epsilon_{0} < \min \left\{\frac{\rho_{0}^{2\nu}}{4^{2\nu}(3c^{3}+5c^{2}+7c+8)^{2}},\frac{1}{4^{2\nu}} \right\}.
\end{equation}
Therefore, conditions \eqref{con11}, \eqref{con22}, and \eqref{cond11} together give the choice of $\epsilon_{0}$.

\noindent\textit{Convergence of $\xi_{j},Z_{j}$:}
Since $\xi$ will be the limit of the sequence $\{\xi_{j}(0)\}_{j}$, it follows that $C(\xi)=\xi$ and $\|\xi\|\leq 2\epsilon_{0}$. Thus, let us show the convergence of $\{\xi_{j}(0)\}_{j}$. By their definition,
$$
\|\xi_{j+1}(0)-\xi_{j}(0)\|=\|\xi_{j}(\xi^{j+1}(0))-\xi_{j}(0)\|\leq \|D\xi_{j}\|_{\sigma_{j+1}}\|\xi^{j+1}(0)\|,
$$
and by the chain rule, for $\|\eta_{j+1}\|\leq \sigma_{j+1}$,
$$
\begin{aligned}\|D\xi_{j}\|_{\sigma_{j+1}}&\leq \|D\xi^{0}\|_{\sigma_{1}}\|D\xi^{1}\|_{\sigma_{2}}\cdots \|D\xi^{j}\|_{\sigma_{j+1}}\\
&\leq (1+c_{1}\sigma_{0}^{\frac{1}{2}})\cdots(1+c_{1}\sigma_{j+1}^{\frac{1}{2}})
\leq \Pi_{j=0}^{\infty}(1+c_{1}\sigma_{j}^{\frac{1}{2}}) \leq \exp \left(c_{1}\sum_{j \geq 0} \sigma_{j}^{\frac{1}{2}} \right).
\end{aligned}
$$
Since $\sigma_{j+1}/\sigma_{j}\leq K\leq \frac{1}{4}$, we have
$$
\sum_{j}\sigma_{j}^{\frac{1}{2}}\leq \sum_{j}\frac{\sigma_{0}^{\frac{1}{2}}}{2^{j}}\leq 2\sigma_{0}^{\frac{1}{2}},
$$
and this implies that $\|D\xi^{j}\|_{\sigma_{j+1}}\leq \infty.$ Therefore, $(\xi_{j}(0))_{j}$ is a Cauchy sequence and it has a limit $\xi$, as desired.

Since $Z_{j}(x,\eta_{j+1})=Z_{j-1}(x,\xi^{j}(\eta_{j+1}))Z^{j}(x,\xi^{j}(\eta_{j+1})),(x,\eta_{j+1})\in D_{j+1}$ for $j\geq 1$ and $Z_{0}(x,\eta_{1})=Z^{0}(x,\xi^{0}(\eta_{1}))$, for $|\im x|<\frac{\rho_{0}}{2}$, we have
$$
\begin{aligned}
&\|Z_{j+1}(x,0)-Z_{j}(x,0)\|\\&=\|Z_{j}(x,\xi^{j+1}(0))Z^{j+1}(x,\xi^{j+1}(0))-Z_{j}(x,0)\|\\
&= \|Z_{j}(x,\xi^{j+1}(0))-Z_{j}(x,0)+Z_{j}(x,\xi^{j+1}(0))(Z^{j+1}(x,\xi^{j+1}(0))-Id)\|\\
&\leq \|Z_{j}(x,\xi^{j+1}(0))-Z_{j}(x,0)\|+\|Z_{j}(x,\xi^{j+1}(0))\|\cdot\|Z^{j+1}(x,\xi^{j+1}(0))-Id\|.
\end{aligned}
$$
Note that $Z^{j+1}=\exp(\chi X_{j+1})$. Thus we have
\begin{eqnarray*}
\|Z^{j+1}(x,\xi^{j+1}(0))-Id\| &\leq& 2\chi\|X_{j+1}\|_{D^{j+2}} \leq 2|\chi|\epsilon_{j+1}^{\frac{1}{2}},\\
\|Z_{j}(x,\xi^{j+1}(0))\|&\leq& \Pi_{j\geq 0}\exp(|\chi| \epsilon_{j}^{\frac{1}{2}}) \leq \exp \left(2|\chi|\epsilon_{0}^{\frac{1}{2}} \right),\\
\|Z_{j}(x,\xi^{j+1}(0))-Z_{j}(x,0)\|&\leq& \|D_{\eta_{j+1}}Z_{j}\|_{D^{j+1}}\|\xi^{j+1}(0)\|.
\end{eqnarray*}
We need the following estimate for $D_{\eta_{j+1}}Z_{j}$. Again by the Cauchy inequality, we have
$$
\begin{aligned}\|D_{\eta_{j+1}}Z_{j}\|_{D^{j+1}}&=\|D_{\eta_{j+1}}Z_{j-1}(\cdot,\xi^{j}(\cdot))Z^{j}(\cdot,\xi^{j}(\cdot))\|_{D^{j+1}}\\
&\leq \|D_{\eta_{j}}(Z_{j-1}Z^{j})\|_{\{|\im x|<\rho_{j+1}\}\times \{|\eta_{j}|<2\sigma_{j+1}\}}\|D\xi^{j}\|_{\sigma_{j+1}}\\
&\leq \frac{c'\|Z_{j-1}Z^{j}\|_{D^{j+1}}}{\sigma_{j}/2-2\sigma_{j+1}}\|D\xi^{j}\|_{\sigma_{j+1}}\\
&\leq \frac{2c'}{\sigma_{j}-4\sigma_{j+1}} \exp \left(2|\chi|\epsilon_{0}^{\frac{1}{2}} \right) \exp \left( |\chi|\sigma_{j+1}^{\frac{1}{2}} \right)(1+c_{1}\sigma_{j})\\
&<\frac{c_{2}}{\sigma_{j}-4\sigma_{j+1}},
\end{aligned}
$$
where $c_{2}$ is a new constant. Therefore, we have
$$\|Z_{j}(x,\xi^{r+1}(0))-Z_{j}(x,0)\|\leq \frac{c_{2}\sigma_{j+1}}{\sigma_{j}-4\sigma_{j+1}}=\frac{c_{2}\sigma_{j}^{\frac{3}{2}}}{\sigma_{j}-4\sigma_{j}^{\frac{3}{2}}}\leq c_{3}\sigma_{j}^{\frac{1}{2}},$$
where $c_{3}$ is again a new constant. Putting all these together, we have
$$
\|Z_{j+1}(x,0)-Z_{j}(x,0)\| \leq \left(c_{3}+2| \chi | \exp \left(2|\chi|\epsilon_{0}^{\frac{1}{2}} \right) \right)\sigma_{j}^{\frac{1}{2}},$$
which implies that $(Z_{j}(x,0))_{j}$ is a Cauchy sequence.

\medskip

\noindent\textit{Analytic dependence on $\mu$:} First of all, the proof above can be extended to $P:\mathbb{T}^{d}\to g_{\mathbb{C}}$ with $\|P\|_{\rho_{0}}\leq \epsilon_{0}$ and to complex $\chi$ with $|\chi|\leq 1$, where $g_{\mathbb{C}}=\{P_{1}+iP_{2}: P_{1}, P_{2}\in g\}$. The admissibility of $(A_{0},C,S,\omega)$ guarantees that for such $P$, there exists a unique analytic solution $X:2\mathbb{T}^{d}\to g_{\mathbb{C}}$ to
$$e^{-A_{0}}[e^{A_{0}},X]=-(P-C([P])),[X]=S([P])$$
with an analytic extension to $|\im x|<\rho_{0}$ and the estimate $\|X\|_{\rho_{0}-\delta}\leq c\frac{\|P\|_{\rho_{0}}}{\delta^{\nu}}$ for all $0<\delta<\rho_{0}$.

Now assume that $P,\chi$ depends analytically on $\mu$ in a neighborhood of zero. Let $\nu>0$ be such that $|\mu|<\nu$ implies $\|P(\cdot,\mu)\|_{\rho_{0}}<\epsilon_{0}$ and $|\chi(\mu)|<1$. For these complex values of $\mu$, there exists $\xi^{*}(\mu)\in \sl(2,\mathbb{R})$ (with $C(\xi^{*})=\xi^{*}$ and $\xi^{*}(\mu)\in \sl(2,\mathbb{R})$ for real values of $\mu$) and $X(\cdot,\mu): 2\mathbb{T}^{d}\to \sl(2,\mathbb{R})$ with analytic extension to $|\im x|<\frac{\rho_{0}}{2}$ such that $Z(x,\mu)=\exp (\chi(\mu)X(x,\mu))$ satisfies
$$
Z^{-1}(x+\omega,\mu)A_{0}e^{\chi^{p+1}P(x,\mu)}e^{-\chi^{p+1}\xi^{*}(\mu)}Z(x+\omega,\mu)=A_{0}
$$
for $|\im x|<\frac{\rho_{0}}{2}$ and $|\mu|<\nu$. Moreover, if $\mu$ is real, then $P\in \sl(2,\mathbb{R})$ and $\xi^{*}(\mu)\in \sl(2,\mathbb{R})$ and $X(x,\mu)\in \sl(2,\mathbb{R})$ for real $x\in 2\mathbb{T}^{d}$. Since the transformations constructed above can be made analytic in $D^{j}\times \{|\mu|<\nu\}$ and the convergence is uniform in the complex domain $D^{*}\times\{|\mu|<\nu\}$, the limits are also analytic there.
\end{pf}

\medskip

\noindent\textbf{Proof of Theorem \ref{KAM1}:}
The proof of this result follows the same lines as that of Theorem~\ref{KAM}, although they have different linearized equations.
That is, there exist $X(x)\in C^{\omega}(\mathbb{T}^{d},g)$ and $Z(x)=\exp(\chi X(x))$ such that
$$
Z^{-1}(x+\omega)\exp(\chi^{p}A_{0}+\chi^{p+1} P(x,\mu)-\chi^{p+1}\xi^{*}(\mu))Z(x)=e^{\chi^{p}A_{0}},
$$
and the linear version of this equation is
$$
-e^{-\chi^{p}A_{0}}X(x+\omega)e^{\chi^{p}A_{0}}+X(x)=-\chi^{p}(\Lambda(P)-C(\Lambda([Q]))),
$$
where $\Lambda(X)=\frac{1}{2}X+B^{-1}XB$ and $B=\begin{pmatrix}1 & 1\\ 0 &1\end{pmatrix}$. Let $X=\begin{pmatrix}X_{11} &X_{12}\\X_{21} &X_{22}\end{pmatrix}$. Then, $\Lambda(X)=\begin{pmatrix}X_{11}-X_{21}/2 &X_{11}/2+X_{12}-X_{21}/2-X_{22}/2\\ X_{21} &X_{22}+X_{21}/2\end{pmatrix}$. Moreover, we have $\Lambda C=C\Lambda$ if $(e^{\chi^{p}A_{0}},C,S,\omega)$ is admissible and $A_{0}=\begin{pmatrix}0 &1\\0 &0\end{pmatrix}$. Assume that
$$
\begin{aligned}&(Z^{j}(x+\omega))^{-1}\exp \chi^{p}(A_{0}+\chi P^{j}(x,\eta^{j})-\chi\eta^{j})Z^{j}(x)\\&=\exp \chi^{p}(A_{0}+\chi P^{j+1}(x,\eta^{j+1})-\chi\eta^{j+1}).\end{aligned}
$$
Then we have
$$
\Lambda(\eta^{j+1})=\Lambda(\eta^{j})-C(\Lambda([P^{j}(x,\eta^{j})])).
$$
Since $\Lambda^{-1}(0)=0$, and $\Lambda, C$ commute, we have
$$
\eta^{j+1}=\eta^{j}-C([P^{j}(x,\eta^{j})]).
$$
This equation is exactly the same as \eqref{eta}, and therefore the remainder of the proof follows from the same arguments as above.

\section{Proof of Lemma~\ref{HYLem}}\label{AB}

We need the following result from \cite{BB}.

\begin{Theorem}\label{IMT}
Let $X,Y,Z$ be Banach spaces, denote by $\mathcal{O}$ the set from which the external parameter $\mu$ is chosen, and suppose that $U\subset X$ and $V \subset Y$ are neighborhoods of $x_{0}\in X$ and $y_{0}\in Y$, respectively. For some $s,t>0$, define $B_{s}(x_{0})=\{x\in X:\|x-x_{0}\|_{X}<s\}, B_{t}(y_{0})=\{y\in Y:\|y-y_{0}\|_{Y}<t\}$. Let $\Phi(x,y;\mu)\in C^{0}(U\times V\times \mathcal{O},Z)$ and for each $\mu\in \mathcal{O},\Phi(\cdot,\cdot;\mu)\in C^{1}(U\times V,Z)$, let $B_{s}(x_{0})\times B_{t}(y_{0})\subset U\times V$, and assume that $\Phi(x_{0},y_{0};\mu)=0$ and $D_{y}\Phi(x_{0},y_{0};\mu)\in \mathcal{L}(Y,Z)$ is invertible. If
$$\sup_{\overline{B_{s}(x_{0})}}\|\Phi(x,y_{0};\mu)\|_{Z}\leq\frac{t}{2\|(D_{y}\Phi(x_{0},y_{0};\mu))^{-1}\|_{Z}},$$
$$\sup_{\overline{B_{s}(x_{0})}\times\overline{B_{t}(y_{0})}} \|Id_{Y}-(D_{y}\Phi(x_{0},y_{0};\mu))^{-1}D_{y}\Phi(x,y;\mu)\|_{\mathcal{L}(Y,Y)}\leq \frac{1}{2},$$
then there exists $y(\cdot,\mu)\in C^{1}(B_{s}(x_{0}),\overline{B_{t}(y_{0})})$ such that $\Phi(x,y(x,\mu);\mu)=0.$ In particular, by Proposition 1.2 of \cite{Zeidler}, we have $y(\cdot,\mu)\in C^{0}(B_{s}(x_{0})\times\mathcal{O},\overline{B_{t}(y_{0})})$.
\end{Theorem}

For $\mu\in \mathcal{O}$, define the non-linear functional
$$
\Phi_{\mu}: \mathcal{B}_{\mathcal{O},r}^{nre}(\mathrm{su}(1,1),\eta)\times \mathcal{B}_{\mathcal{O},r}(\mathrm{su}(1,1))\to \mathcal{B}_{\mathcal{O},r}^{nre}(\mathrm{su}(1,1),\eta)$$
by $$\Phi_{\mu}(Y,g)=\mathbb{P}_{nre}\ln(e^{A^{-1}(\mu)Y(\mu,x+\tilde{\omega})A(\mu)}e^{g(\mu,x)}e^{-Y(\mu,x)}).
$$
Clearly, we have $\Phi_{\mu}(0,0)=0$ and $\|\Phi_{\mu}(0,g)\|\leq \|g\|_{r}$. Now we consider the difference
$$
\begin{aligned}&\Phi_{\mu}(Y+Y',g)-\Phi_{\mu}(Y,g)\\
&=\mathbb{P}_{nre}\ln(e^{A^{-1}(\mu)(Y(\mu,x+\tilde{\omega})+Y'(\mu,x+\tilde{\omega}))A(\mu)}e^{g(\mu,x)}e^{-Y(\mu,x)-Y'(\mu,x)})\\
&-\mathbb{P}_{nre}\ln(e^{A^{-1}(\mu)Y(\mu,x+\tilde{\omega})A(\mu)}e^{g(\mu,x)}e^{-Y(\mu,x)})\\
&=\mathbb{P}_{nre}\ln(e^{A^{-1}(\mu)(Y(\mu,x+\tilde{\omega})+Y'(\mu,x+\tilde{\omega}))A(\mu)}e^{g(\mu,x)}e^{-Y(\mu,x)-Y'(\mu,x)})\\
&-\mathbb{P}_{nre}\ln(e^{A^{-1}(\mu)(Y(\mu,x+\tilde{\omega})+Y'(\mu,x+\tilde{\omega}))A(\mu)}e^{g(\mu,x)}e^{-Y(\mu,x)})\\
&+\mathbb{P}_{nre}\ln(e^{A^{-1}(\mu)(Y(\mu,x+\tilde{\omega})+Y'(\mu,x+\tilde{\omega}))A(\mu)}e^{g(\mu,x)}e^{-Y(\mu,x)})\\
&-\mathbb{P}_{nre}\ln(e^{A^{-1}(\mu)Y(\mu,x+\tilde{\omega})A(\mu)}e^{g(\mu,x)}e^{-Y(\mu,x)})\\
\end{aligned}
$$
We need the Baker-Campbell-Hausdorff formula:
$$
\ln(e^{X}e^{Y})=X+Y+\frac{1}{2}[X,Y]+\frac{1}{12}\{[X,[X,Y]]+[Y,[X,Y]]\}+\cdots,
$$
where $[X,Y]=XY-YX$ denotes the Lie bracket and $\cdots$ denotes the higher order terms of Lie brackets. Then we have
$$
\begin{aligned}
&\Phi_{\mu}(Y+Y',g)-\Phi_{\mu}(Y,g)\\&=\mathbb{P}_{nre}(-Y'+A^{-1}(\mu)Y'(\mu,x+\tilde{\omega})A(\mu)+O_{2}(Y,Y',g)).
\end
{aligned}
$$
Let $Y=0$ and $g=0$. Then the higher order terms vanish and we obtain
$$
\begin{aligned}
D_{Y}\Phi_{\mu}(0,0)(Y')&=\mathbb{P}_{nre}(-Y'(\mu,x)+A^{-1}(\mu)Y'(\mu,x+\tilde{\omega})A(\mu))\\
&=-Y'(\mu,x)+A^{-1}(\mu)Y'(\mu,x+\tilde{\omega})A(\mu),
\end{aligned}
$$
Therefore, we have
$$
\|D_{Y}\Phi_{\mu}(0,0)(Y')\|\geq\|-Y'(\mu,x)+A^{-1}(\mu)Y(\mu,x+\tilde{\omega})A(\mu)\|_{r}\geq \eta\|Y\|_{r}
$$
and $\|D_{Y}\Phi_{\mu}(0,0)^{-1}\|\leq \eta^{-1}$.

Choose $s=\epsilon$, $t=\epsilon^{\frac{1}{2}}$, and $\eta\geq 13\|A_{\mu}\|^{2}\epsilon^{\frac{1}{2}}$. We have
$$
2\sup_{\overline{B_{s}(x_{0})}}\|\Phi_{\mu}(0,g)\|\times\|D_{Y}\Phi_{\mu}(0,0)^{-1}\|\leq 2\epsilon\times\frac{\epsilon^{-\frac{1}{2}}}{13\|A(\mu)\|^{2}}\leq\epsilon^{\frac{1}{2}},
$$
and this verifies the first condition of Theorem \ref{IMT}. In order to verify the second condition, note that
$$
\begin{aligned}
\|D_{Y}\Phi_{\mu}(Y,g)(Y')-D_{Y}\Phi_{\mu}(0,0)(Y')\|&=\|O(Y,g)Y'\|_{r}\\
&\leq 6(\|A(\mu)\|^{2}\|Y\|_{r}+\|A(\mu)\|^{2}\|g\|_{r})\|Y'\|_{r}\\
&\leq 6\|A(\mu)\|^{2}(\epsilon^{\frac{1}{2}}+\epsilon)\|Y'\|_{r}\end{aligned}.
$$
This implies that
$$
\sup_{B_{s}(x_{0})\times B_{t}(y_{0})}\|D_{Y}\Phi_{\mu}(0,0)-D_{Y}\Phi_{\mu}(Y,g)\|\leq 6\|A(\mu)\|^{2}(\epsilon^{\frac{1}{2}}+\epsilon),
$$
and therefore
$$
\sup_{B_{s}(x_{0})\times B_{t}(y_{0})}\|Id_{\mathcal{B}_{nre}(\eta,\mu)}-D_{Y}\Phi_{\mu}(0,0)^{-1}D_{Y}\Phi_{\mu}(Y,g)\|\leq \frac{1}{2},
$$
and the second condition holds as well.

Now we can apply Theorem \ref{IMT} to $\|g\|_{r}\leq\epsilon, \eta\geq 13\|A(\mu)\|^{2}\epsilon^{\frac{1}{2}}$ and obtain $Y(\mu,\cdot)$ with $\|Y(\mu,\cdot)\|_{r}\leq \epsilon^{\frac{1}{2}}$, which satisfies
$\Phi_{\mu}(Y(\mu,\cdot),g(\mu,\cdot))=0$, that is,
$$
e^{A^{-1}(\mu)Y(\mu,x+\tilde{\omega})A(\mu)}e^{g(\mu,x)}e^{-Y(\mu,x)}=e^{g^{re}(\mu,x)},
$$
which is equivalent to
$$
e^{Y(\mu,x+\tilde{\omega})}A(\mu)e^{g(\mu,x)}e^{-Y(\mu,x)}=A(\mu)e^{g^{re}(\mu,x)}.
$$

\end{document}